%% file: sym2p2-2008.07.24.tex
\documentclass[10pt,letterpaper]{article}

\usepackage{fullpage}
\usepackage{defs}
\usepackage{words}
\usepackage{art-thm}
\usepackage{url}

\renewcommand{\labelenumi}{(\alph{enumi})}

\numberwithin{equation}{subsection}

\providecommand{\GW}[2]{\ensuremath{\langle {#1} \rangle_{{#2}}}}

\begin{document}

\title{The genus zero Gromov--Witten invariants of $[\Sym^2 {\bf P}^2]$}
\author{Jonathan Wise}
\date{\today}

\maketitle

\begin{abstract}
  We study the Abramovich--Vistoli moduli space of genus zero orbifold stable maps to $[\Sym^2 \PP^2]$, the stack symmetric square of $\PP^2$.  This compactifies the moduli space of stable maps from hyperelliptic curves to $\PP^2$, and we show that all genus zero Gromov--Witten invariants are determined from trivial enumerative geometry of hyperelliptic curves.  We also show how the genus zero Gromov--Witten invariants can be used to determine the number of hyperelliptic curves of degree $d$ and genus $g$ interpolating $3d + 1$ generic points in $\PP^2$.  Comparing our method to that of Graber for calculating the same numbers, we verify an example of the crepant resolution conjecture.
\end{abstract}

\tableofcontents

\newcommand{\orbifoldRRref}{\cite{AGV}, Section~7.2}
\newcommand{\reconstructionproofref}{Appendix~\ref{reconstruction-proof}}
\renewcommand{\theparagraph}{\arabic{section}.\arabic{subsection}.\arabic{paragraph}}
\newcommand{\hodgeintegralreference}{\cite{hhi}}

\section{Introduction}

\input{intro-2008.07.21.tex}

\section{Basic facts about $[\Sym^2 \PP^2]$}

\input{basic-2008.07.21.tex}

\section{Moduli of orbifold stable maps}
\label{moduli}
\providecommand{\modulisection}{\subsection}
\input{moduli-2008.07.24.tex}

\section{Gromov--Witten invariants and enumerative geometry}
\label{sect:GW-enum}
\label{gw-enum}
\input{gw-enum-2008.07.24.tex}

\section{Calculating the Gromov--Witten invariants}
\label{gw}
\input{gw-2008.06.15.tex}

\section{The crepant resolution conjecture}
\label{sect:crc}
\label{crc}
\input{crc-2008.06.15.tex}

\bibliographystyle{plain}
\bibliography{thesis}

\end{document}

%% file: intro-2008.07.21.tex
Abramovich and Vistoli were motivated in their definition of the moduli space of orbifold stable maps to compactify the space of stable maps to an orbifold.  The original definition of orbifold Gromov--Witten theory, by Chen and Ruan~\cite{CR}, was motivated by mirror symmetry in dimensions larger than $3$.  However, neither of these motivations comes into play here.  We are motivated instead by Graber's enumeration of hyperelliptic curves in $\PP^2$~\cite{G}: by viewing hyperelliptic curves in $\PP^2$ as families of length $2$ subschemes of $\PP^2$ parameterized by a rational curve, the space of stable maps to $\Hilb_2 \PP^2$ becomes a compactification of the moduli space of hyperelliptic curves in $\PP^2$.  

Now that the Abramovich--Vistoli moduli space is available, an even more natural compactification presents itself.  A hyperelliptic curve in $\PP^2$ is nothing but a family of pairs of points in $\PP^2$, parameterized by an orbifold curve of genus zero.  Thus the genus zero Gromov--Witten invariants of $[\Sym^2 \PP^2]$ capture the enumerative geometry of hyperelliptic curves in $\PP^2$ with the only twist being the presence of the virtual fundamental class.


By definition,
\begin{equation*}
  [\Sym^2 {\bf P}^2] = \left[ \left. \left( {\bf P}^2 \times {\bf P}^2 \right) \right/ S_2   \right],
\end{equation*}
with $S_2$ acting by exchanging the components.  This is a smooth, $4$-dimensional Deligne--Mumford stack whose coarse moduli space $\Sym^2 {\bf P}^2$ has an $A_1$ surface singularity along the diagonal.  Note that we will always work over $\CC$ in this paper, and therefore we identify $S_2 \cong \ZZ / 2 \ZZ \cong \mu_2$ without further comment.  

The Abramovich--Vistoli moduli space associated to a Deligne--Mumford stack $X$ will be denoted here by $\overline{M}(X)$ with various decorations to specify connected components.  A point of $\overline{M}(X)$ corresponds to a representable morphism $C \rightarrow X$ (with a number of other properties).  When $X = [\Sym^2 \PP^2]$, this means that if $\tilde{C}$ is defined to make the diagram
\begin{equation*} \xymatrix{
    \tilde{C} \ar[r]^<>(0.5){\tilde{f}} \ar[d] & \PP^2 \times \PP^2 \ar[d] \\
    C \ar[r]^<>(0.5)f & [\Sym^2 \PP^2]
} \end{equation*}
cartesian, then $\tilde{C}$ is a scheme posessing an $S_2$ action with respect to which $\tilde{f}$ is equivariant.  The coase quotient of $\tilde{C}$ by this $S_2$-action is the coarse moduli space of $C$;  if $C$ is a smooth curve of genus zero (meaning its coarse moduli space has genus zero), this implies that $\tilde{C}$ is a hyperelliptic curve.  Moreover, the equivariant map $\tilde{f}$ is determined by its projection on either factor.  Conversely, any map from a hyperelliptic curve to $\PP^2$ induces an equivariant map to $\PP^2 \times \PP^2$, and thus we see that there is an open substack of $\overline{M}([\Sym^2 \PP^2])$ parameterizing stable maps from hyperelliptic curves to $\PP^2$.  In other words, we have demonstrated that $\overline{M}([\Sym^2 \PP^2])$ is indeed a compactification of the moduli space of stable hyperelliptic curves in $\PP^2$.  We shall therefore find the enumerative geometry of hyperelliptic curves in $\PP^2$ reflected in the genus zero Gromov--Witten theory of $[\Sym^2 \PP^2]$.

We will begin our study of the moduli space $\overline{M}([\Sym^2 \PP^2])$ in Section~\ref{moduli}.  We will not achieve a complete description of the moduli space in any sense, but we will at least identify a collection of irreducible components that are sufficient to address the enumerative problems of Section~\ref{gw-enum}.  The main results of Section~\ref{moduli} are Theorem~\ref{dim-ev-image}, which implies that for the purpose of counting hyperelliptic curves through points in $\PP^2$, one may restrict one's attention to the open locus of comb curves in $\overline{M}([\Sym^2 \PP^2])$, and Theorem~\ref{thm:mod:1}, which evaluates the contribution of each of these components.

This will permit us, in Section~\ref{gw-enum}, to reduce the enumeration of hyperelliptic curves passing through point in $\PP^2$ to one of counting the connected components of the moduli space of such curves.  We prove the following there as Theorem~\ref{thm:2}.

\begin{mainthm} \label{mainthm:2}
  Let $E(d,g)$ be the number of hyperelliptic curves in $\PP^2$ passing through $3d + 1$ points in generic position and let $J(d,g)$ be the corresponding Gromov--Witten invariant of $[\Sym^2 \PP^2]$.  Then
  \begin{equation*}
    J(d,g) = \sum_{{\bf h}} \left( - \frac{1}{4} \right)^{g - g({\bf h})} (2 g({{\bf h}}) + 2)! E(d, g({\bf h})) .
  \end{equation*}
  The sum is taken all partitions of $[2g + 2]$ into $2 g({\bf h}) + 2$ parts having odd numbers of elements.
\end{mainthm}
The factor of $\left( - \frac{1}{4} \right)^{g - g({\bf h})}$ is explained by a hyperelliptic Hodge integral.  It seemed too distracting to include this calculation here, so it will appear elsewhere~\cite{hhi}.

Section~\ref{gw-enum} also gives the relationships between several other genus zero Gromov--Witten invariants of $[\Sym^2 \PP^2]$ and corresponding enumerative problems.  In Section~\ref{gw} we prove that these are enough to determine all of the genus zero Gromov--Witten invariants of $[\Sym^2 \PP^2]$ by means of the WDVV equations, in a manner essentially the same as the proof of the Kontsevich--Manin reconstruction theorem.  As the expression in Thereom~\ref{mainthm:2} can be inverted to express the $E(d,g)$ in terms of the $J(d,g)$, we obtain an algorithm to determine the $E(d,g)$ recursively.

Finally, in Section~\ref{crc}, we compare our methods to those of Graber in the manner suggested by Ruan's crepant resolution conjecture~\cite{R},~\cite{BG},~\cite{CoR}.  Of course, both approaches give the same answer to the enumerative problem in the end.  However, each approach requires the evaluation of Gromov--Witten invariants with non-trivial contributions from the two different compactifications.  Ruan's conjecture predicts that there should be a direct relationship on the level of Gromov--Witten theory.  We verify this in Section~\ref{crc}.

\begin{mainthm} \label{mainthm:3}
  The crepant resolution conjecture is valid as stated by Bryan and Graber~\cite{BG} for the resolutions $[\Sym^2 \PP^2]$ and $\Hilb_2 \PP^2$ of $\Sym^2 \PP^2$.
\end{mainthm}

\subsection{Acknowledgements}

My work on this project has benefitted from conversations with Jim Bryan, Charles Cadman, Barbara Fantechi, W. D. Gillam, Tom Graber, Rahul Pandharipande, Angelo Vistoli, and Ben Wieland.  I am especially grateful to my advisor, Dan Abramovich, for his unending patience and encouragement.  

I also thank the Institut Henri Poinar\'e for its hospitality while a portion of this paper was being written.

%% file: basic-2008.07.21.tex
Here we gather some properties of $[\Sym^2 \PP^2]$ that do not specifically concern curves, but which we will need later.

The stack $[\Sym^2 \PP^2]$ is the moduli space of unordered pairs of points on $\PP^2$.  More precisely, it is the stack quotient, $[ \PP^2 \times \PP^2 / S_2 ]$ with $S_2$ acting by exchanging the factors.  We will frequently use $\pi$ to denote the $2$-to-$1$ \'etale cover $\PP^2 \times \PP^2 \rightarrow [\Sym^2 \PP^2]$.  There is also a canonical map $[\Sym^2 \PP^2] \rightarrow B S_2$ induced by the equivariant map from $\PP^2 \times \PP^2$ to a point.

The equivariant embedding of the diagonal (with the trivial action of $S_2$) in $\PP^2 \times \PP^2$ induces a closed substack $\Delta \subset [\Sym^2 \PP^2]$, which we also call the diagonal.  It is isomorphic to $\PP^2 \times B S_2$.

Typically, if $Z \subset [\Sym^2 \PP^2]$ is a substack, we will write $\tilde{Z} = \pi^{-1}(Z)$ for its pullback to $\PP^2 \times \PP^2$.  Thus $\tilde{\Delta}$ is the diagonal in $\PP^2 \times \PP^2$.  More generally, though less precisely, we will apply a tilde to a construction for $[\Sym^2 \PP^2]$ to denote a corresponding construction for $\PP^2 \times \PP^2$.  To describe a $P$-point of $[\Sym^2 \PP^2]$, we will often give $\tilde{P}$ with an $S_2$ action and an equivariant map $\tilde{P} \rightarrow \PP^2 \times \PP^2$.  

Note finally that the action of $\PGL{3}$ on $\PP^2$ is $2$-transitive, so the induced action on $[\Sym^2 \PP^2]$ has two orbits: the diagonal and its complement.

\subsection{Standard vector bundles}
\label{vector-bundles}

We define vector bundles $E_1$ and $E_2$ on $[\Sym^2 \PP^2]$ by writing their pullback to $\PP^2 \times \PP^2$ and giving the induced action of $S_2$.  The vector bundles are
\begin{equation*} \begin{split}
    \pi^\ast(E_1) & = {\cal O}(1) \boxtimes {\cal O}(1) = p_1^\ast {\cal O}_{\PP^2}(1) \tensor p_2^\ast {\cal O}_{\PP^2}(1) \\
    \pi^\ast(E_2) & = {\cal O}(1) \boxplus {\cal O}(1) = p_1^\ast {\cal O}_{\PP^2}(1) \oplus p_2^\ast {\cal O}_{\PP^2}(1)
  \end{split}
\end{equation*}
with $S_2$ acting on each by exchaging the components.  We also have line bundles $\rho_0$ and $\rho_1$, the trivial and non-trivial representations of $S_2$, respectively, pulled back from the canonical morphism $[\Sym^2 \PP^2] \rightarrow B S_2$.  Of course, $\rho_0$ is the trivial line bundle on $[\Sym^2 \PP^2]$ and we also denote it by ${\cal O}$ (making the standard identification between line bundles and invertible sheaves).

A global section of $E_1$ may be viewed as a polynomial of bihomogeneous degree $(1,1)$ in two sets of $3$ variables that is invariant under the exchange of the two sets of variables.  If $x$ is a coordinate on $\PP^2$ vanishing along a hyperplane $H$ then $x \tensor x$ determines a section of $E_1$ that vanishes on $[ (H \times \PP^2 \cup \PP^2 \times H) / S_2 ] \subset [ \Sym^2 \PP^2]$.

A global section of $\pi^\ast(E_2)$ is a pair of sections of ${\cal O}_{\PP^2}(1)$.  Sections of $E_2$ are those pairs that are invariant under the action of $S_2$; these can be identified with sections of ${\cal O}_{\PP^2}(1)$ over $\PP^2$.  In general, the vanishing locus of a section of $E_2$ is $[ \Sym^2 H ] \cong [ \Sym^2 \PP^1 ]$ where $H \cong \PP^1$ is a line in $\PP^2$.  As it will be important later, we note that $[ \Sym^2 H ]$ has codimension $2$ in $[\Sym^2 \PP^2]$ but intersects the diagonal in codimension $1$.

On a stack, a vector bundle whose fiber at a stacky point has a non-trivial action of the stabilizer group cannot be generated by global sections, since stabilizers act trivially on global sections.  The most we can hope is that $H^0(P, F \rest{P})$ should be generated by global sections.  This is the case for $E_2$.

\begin{proposition}
  If $P$ is a zero dimensional integral closed substack of $[\Sym^2 \PP^2]$, then 
\begin{equation*}
H^0([\Sym^2 \PP^2], E_2) \rightarrow H^0([\Sym^2 \PP^2], E_2 \rest{P})
\end{equation*}
is surjective.
\end{proposition}
\begin{proof}
  We prove the proposition for $P \in \Delta$ and $P \not\in \Delta$ separately.  If $P \subset \Delta$ is represented by $(p, p) \in \tilde{\Delta}$, then $\dim H^0(P, E_2 \rest{P}) = 1$.  Identifying sections of $E_2$ over $[\Sym^2 \PP^2]$ with sections of ${\cal O}(1)$ over $\PP^2$, it is thus sufficient to find a section of ${\cal O}(1)$ that does not vanish at $p$.

  If $P$ is represented by $(p,q) \amalg (q,p) \not\in \tilde{\Delta}$, then
\begin{equation*}
H^0(P, E_2 \rest{P}) = H^0(\tilde{P}, \pi^\ast E_2)^{S_2} \cong H^0(p, {\cal O}_{\PP^2}(1) \rest{p}) \oplus H^0(q, {\cal O}_{\PP^2}(1) \rest{q}).
\end{equation*}
We may certainly find a pair of sections of ${\cal O}_{\PP^2}(1)$ such that one vanishes at $p$ but not at $q$ and the other vanishes at $q$ but not at $p$, so the proof is complete.
\end{proof}

\begin{corollary} \label{glob-gen-of-tan}
If $P$ is an integral closed substack of dimension zero in $[\Sym^2 \PP^2]$ then 
\begin{equation*}
	H^0([\Sym^2 \PP^2], T [\Sym^2 \PP^2]) \rightarrow H^0(P, T [\Sym^2 \PP^2] \rest{P})
\end{equation*}
is surjective.
\end{corollary}
\begin{proof}
The Euler sequence on $\PP^2$ is
\begin{equation*}
  0 \rightarrow {\cal O} \rightarrow {\cal O}(1)^{\oplus 3} \rightarrow T \PP^2 \rightarrow 0 .
\end{equation*}
This induces an exact sequence on $[ \Sym^2 \PP^2]$,
\begin{equation*}
  0 \rightarrow \rho_0 \oplus \rho_1 \rightarrow E_2^{\oplus 3} \rightarrow T [\Sym^2 \PP^2] \rightarrow 0 .
\end{equation*}
We have the commutative diagram
\begin{equation*} \xymatrix{
H^0([\Sym^2 \PP^2], E_2^{\oplus 3}) \ar[r] \ar[d] & H^0(P, E_2^{\oplus 3} \rest{P}) \ar[d] \\
H^0([\Sym^2 \PP^2], T [\Sym^2 \PP^2]) \ar[r] & H^0(P, T [\Sym^2 \PP^2] \rest{P})
} \end{equation*}
where the vertical arrows come from the Euler sequences.  The upper horizontal arrow is surjective by the proposition; the vertical arrow on the right is surjective because taking global sections over $P$ corresponds to taking $S_2$-invariants, which is exact in characteristic zero.  Thus the lower horizontal arrow must therefore be surjective as well.
\end{proof}

\subsection{The inertia stack}

By definition, the points of the inertia stack $I [\Sym^2 \PP^2]$ are pairs $(x, g)$ where $x$ is a point of $[\Sym^2 \PP^2]$ and $g$ is an automorphism of $x$.  Therefore $I [\Sym^2 \PP^2]$ has two components,
\begin{equation*} \begin{split}
    I [\Sym^2 \PP^2] & = \Omega_0 \amalg \Omega_1 \\
    \Omega_0 & \cong [\Sym^2 \PP^2] \\
    \Omega_1 & \cong \Delta \cong \PP^2 \times B S_2 .
  \end{split}
\end{equation*}
The inertia stack classifies maps from trivialized gerbes under cyclic groups into $[\Sym^2 \PP^2]$.  

If $I_r [\Sym^2 \PP^2]$ is the component where $g$ has order $r$, then $I_r [\Sym^2 \PP^2]$ has a natural faithful action of $B (\ZZ / r \ZZ)$.  The quotient by this action is called the rigidified $\ZZ / r \ZZ$-inertia stack and is written $\overline{I}_r [\Sym^2 \PP^2]$.  It classifies maps from gerbes banded by $\ZZ / r \ZZ$ into $[\Sym^2 \PP^2]$.  The total rigidified inertia stack is the disjoint union of the $\overline{I}_r [\Sym^2 \PP^2]$:
\begin{equation*} \begin{split}
    \overline{I} [\Sym^2 \PP^2] & = \overline{\Omega}_0 \amalg \overline{\Omega}_1 \\
    \overline{\Omega}_0 & = [\Sym^2 \PP^2] \\
    \overline{\Omega}_1 & = \PP^2 .
    \end{split}
\end{equation*}

\subsection{The orbifold Chow group}
\label{group-structure}

By definition, the orbifold Chow ring is 
\begin{equation*}
  A_{\rm orb}^\ast (\overline{I} [\Sym^2 \PP^2]) = A^\ast ([\Sym^2 \PP^2]) \oplus A^\ast(\PP^2)
\end{equation*}
as a vector space (throughout this paper, all Chow groups will be taken with rational coefficients).  Its grading is shifted by the {\em age} which is discussed below.  

We will determine the orbifold Chow ring of $[\Sym^2 \PP^2]$ as a graded vector space in this section and defer the discussion of its ring structure until we have computed the Gromov--Witten invariants of $[\Sym^2 \PP^2]$ that are needed for its definition (at least in the usual way, but see also \cite{JKK}).

\subsubsection{Group structure}
\label{age}

As $[\Sym^2 \PP^2]$ is smooth, the underlying vector space of its Chow ring may be identified with its Chow group.  The map
\begin{equation*}
\pi^\ast : A^\ast([\Sym^2 {\bf P}^2]) \rightarrow A^\ast({\bf P}^2 \times {\bf P}^2)
\end{equation*}
carries $A^\ast([\Sym^2 {\bf P}^2])$ to the subring of invariants of $A^\ast({\bf P}^2 \times {\bf P}^2)$ under the action of $S_2$ induced by switching the components.  The Chow ring of ${\bf P}^2 \times {\bf P}^2$ is
\begin{equation*}
A^\ast({\bf P}^2 \times {\bf P}^2) = {\bf Q}[h_1, h_2] / (h_1^3, h_2^3),
\end{equation*}
with $h_1 = c_1 \left( p_1^\ast {\cal O}(1) \right)$ and $h_2 = c_1 \left( p_2^\ast {\cal O}(1) \right)$.  It is not difficult to show that the ring of invariants is
\begin{equation*}
\begin{split}
A^\ast([\Sym^2 {\bf P}^2]) 
& \cong {\bf Q}[\alpha, \beta] / (\alpha^3 - 3 \alpha \beta, \alpha^2 \beta - 2 \beta^2, \alpha \beta^2, \beta^3) \\
& = {\bf Q} + {\bf Q} \alpha + {\bf Q} \alpha^2 + {\bf Q}\beta + {\bf Q} \alpha^3 + {\bf Q} \alpha^4
\end{split}
\end{equation*}
with $\pi^\ast(\alpha) = h_1 + h_2$ and $\pi^\ast(\beta) = h_1 h_2$.  

The orbifold Chow group of $[\Sym^2 {\bf P}^2]$ is isomorphic to Chow group of 
\begin{equation*}
\overline{I}(\Sym^2 {\bf P}^2) = \overline{\Omega}_0 \amalg \overline{\Omega}_1.
\end{equation*}
We have $\overline{\Omega}_1 \cong {\bf P}^2$, so its Chow group is
\begin{equation*}
{\bf Q} \gamma_0 + {\bf Q} \gamma_1 + {\bf Q} \gamma_2
\end{equation*}
where $\gamma_0$ is the fundamental class, $\gamma_1$ is the class of a line, and $\gamma_2$ is the class of a point.

\subsubsection{Grading}

If $\gamma$ is a class in $A^p(\overline{\Omega}_i)$ then its {\it orbifold degree} is, by definition, $p + \age(\overline{\Omega}_i)$.  

We recall the definition of the age.  A point of $\overline{I} X$ is a pair $(x, g)$ where $x \in X$ and $g \in \Aut(x)$, determined uniquely up to conjugation by other automorphisms of $x$.  The eigenvalues of the action of $g$ on $T_x X$ are therefore well-defined.  As $g$ has finite order, they are roots of unity, say $e^{2 \pi i t_j}$, $j = 1, \ldots, n$.  Then the age of $(x,g)$ is defined to be $\sum t_j$.  This is a locally constant function on $\overline{I} X$, so we may refer to the age of a component of $\overline{I} X$.

We return to the case ${\scr X} = [\Sym^2 {\bf P}^2]$.  The age of $\overline{\Omega}_0$ is of course zero.  If $(x, g)$ represents a point of $\overline{\Omega}_1$, then we can represent $T_x [\Sym^2 \PP^2]$ as the tangent bundle of $\PP^2 \times \PP^2$ with $g$ acting by exchanging the components.  The eigenvalues are $\pm 1$, each with multiplicity $2$, so $\age(\overline{\Omega}_1) = 1$.

We can now write down $A^\ast_{\rm orb}([\Sym^2 {\bf P}^2])$ as a graded vector space.  It is
\begin{equation*} \begin{array}{cc|c|c|c|c}
    & \vphantom{\Big|} A^0_{\rm orb} & A^1_{\rm orb} & A^2_{\rm orb} & A^3_{\rm orb} & A^4_{\rm orb} \\
    \cline{2-6}
    \vphantom{\Big|} \overline{\Omega}_0 & {\bf Q} & {\bf Q} \alpha & {\bf Q} \alpha^2 + {\bf Q} \beta & {\bf Q} \alpha^3 & {\bf Q} \alpha^4 \\
    \overline{\Omega}_1 & & {\bf Q} \gamma_0 & {\bf Q} \gamma_1 & {\bf Q} \gamma_2
\end{array} \end{equation*}

\subsection{Algebraic equivalence classes of curves}

We have seen in the last section that $A_1([\Sym^2 \PP^2])$ is $1$-dimensional.  Numerical classes of curves may therefore be identified with non-negative integers.  To a curve $C$ of degree $d$ in $[\Sym^2 \PP^2]$ we may associate a corresponding hyperelliptic curve $\tilde{C}$ (recall our convention that $\tilde{C} = C \fprod[{[\Sym^2 \PP^2]}] (\PP^2 \times \PP^2)$) that is equivariantly embedded in $\PP^2 \times \PP^2$.  Composing $\tilde{C} \rightarrow \PP^2 \times \PP^2$ with either of the two projections is thus a curve of degree $d$ in $\PP^2$.

The degree of $f : C \rightarrow [\Sym^2 \PP^2]$ can be defined in a somewhat more intrinsic fashion as the degree of $f^\ast E_1$ where $E_1$ is the line bundle defined in Section~\ref{vector-bundles}.  Indeed,
\begin{equation*}
        \deg f^\ast E_1 = \int_C c_1(f^\ast E_1) = \frac{1}{2} \int_{\tilde{C}} c_1(\tilde{f}^\ast {\cal O}(1,1)) = \frac{1}{2} ( \deg(p_1 \tilde{f}) + \deg(p_2 \tilde{f} ) ) = d
\end{equation*}
where $d = \deg(p_1 \tilde{f}) = \deg(p_2 \tilde{f})$.

%% file: moduli-2008.07.24.tex
\providecommand{\modulisection}{\section}

In this section, we study the moduli space of orbifold stable maps to $[\Sym^2 {\bf P}^2]$.  In Section~\ref{vdim-calcs}, we compute the virtual dimensions of this moduli space.  Sections~\ref{comb-curves} through~\ref{proof-completion} are devoted to the statement and proof of Theorem~\ref{thm:mod:1}, which identifies some of the irreducible components of the moduli space associated to comb curves and determines their virtual fundamental classes.  Section~\ref{evaluation} contains some tedious dimension estimates that legitimize restricting our attention to comb curves in the enumerative applications of Section~\ref{gw-enum}.

\modulisection{Notation}

Let $X$ be a Deligne--Mumford stack with a representable morphism to $B S_2$.  Let $d \in H_2(X, \ZZ)$ be an effective curve class.  We write 
\begin{equation*}
\overline{M}_n(X, d, g)
\end{equation*}
for the moduli space of degree $d$ orbifold stable maps (or twisted stable maps in \cite{AV}) to $X$ with $2g + 2$ orbifold marked points (these must all have automorphism group $S_2$ because we are working over $B S_2$) and $n$ ordinary marked points.  We write $M_n(X, d, g)$ for the open substack of $\overline{M}_n(X, d, g)$ parameterizing orbifold stable maps with smooth source curves.  We also have occasional use for ${\frk M}_n(X, d, g)$, the Artin stack of pre-stable maps to $X$ of degree $d$ with $n$ ordinary marked points and $2g + 2$ orbifold marked points.  If the dimension of $X$ is zero, we omit $d$ from the notation.

We recall that a family of representable maps $C \rightarrow X$ over a base $B$ is called an orbifold pre-stable map if its fibers are nodal orbifold curves (Deligne--Mumford stacks of dimension $1$ with trivial generic stabilizers), with stack structure appearing only at orbifold marked points and at the non-smooth locus of $C \rightarrow B$.  An orbifold marked ``point'' is actually a integral closed substack of $C$ that is a gerbe over $B$ under a cyclic group.  An orbifold pre-stable map is called stable if its automorphism group is finite.

Since our primary concern is $[\Sym^2 \PP^2]$, I will frequently omit $X$ from the notation above when $X = [\Sym^2 \PP^2]$ and there is no danger of confusion.  When $X = [\Sym^2 \PP^2]$, we also use $M^\circ_n(d,g)$ for the open substack of $M_n(d,g)$ consisting of curves that meet the diagonal transversally and only at orbifold points.

\modulisection{Virtual dimension}
\label{sect:vdim-of-sym2p2}
\label{vdim-calcs}

Assume that $X$ is a smooth Deligne--Mumford stack.  It is easiest to define the virtual dimension of $\overline{M}_n(X, d, g)$ relative to ${\frk M}_n(B S_2, g)$ at a point $(C, f)$ corresponding to a map $f : C \rightarrow X$.  In this case, the expected relative dimension is
\begin{equation*}
  \chi(f^\ast T X) = \dim H^0(f^\ast T X) - \dim H^1 (f^\ast T X) .
\end{equation*}
This can be computed by the orbifold Riemann--Roch formula (\orbifoldRRref), which gives
\begin{equation*} \begin{split}
  \chi(f^\ast T X) & = \rank (TX) ( 1 - g(C) ) + \int_C c_1(TX) - \sum_{P \in C} \age_P(f^\ast TX)
\end{split} \end{equation*}
We specialize to the case where $X = [\Sym^2 \PP^2]$ and $g(C) = 0$.  The age of $f^\ast TX$ at an orbifold point is $1$ (as computed in Section~\ref{age}) and $\int_C c_1(TX) = 3d$ so the formula becomes
\begin{equation*}
  \chi(f^\ast T [\Sym^2 \PP^2]) = 4 + 3d - (2g + 2) = 3d - 2g + 2 .
\end{equation*}
The dimension of ${\frk M}_0(B S_2, g)$ is $2g - 1$, so adding this and the contribution of the ordinary marked points gives
\begin{equation*} \begin{split}
  \vdim \overline{M}_{0}([\Sym^2 \PP^2], d, g) & = (3d - 2g + 2) + (2g - 1) + n = 3d + 1 + n \\
\end{split} \end{equation*}

We can also compute this when $X = \Delta$ (which is homogeneous, so in fact the expected relative dimension equals the virtual dimension).  In this case, the age at each orbifold point is now zero since the automorphisms act trivially on the tangent bundle of $\Delta$.  We have $\int_C f^\ast T\Delta = \frac{3d}{2}$ (note that the degree of any map from a curve with trivial generic stabilizer to $\Delta$ must have even degree) and so
\begin{equation*} \begin{split}
    \chi(f^\ast T \Delta) & = 2 + \frac{3d}{2} \\
    \vdim \overline{M}_{n}(\Delta, d, g) & = \frac{3d}{2} + 2g + 1 + n .
  \end{split}
\end{equation*}
Of course, $\vdim \overline{M}_{0}(\Delta, d, g)$ coincides with the dimension of $\overline{M}_{0, 2g + 2}(\PP^2, \frac{d}{2})$ because $\Delta$ is an $S_2$-gerbe over $\PP^2$.

\modulisection{Comb curves} 
\label{comb-curves}

\begin{definition}
A {\em comb curve} is an orbifold stable map $f : C \rightarrow [\Sym^2 \PP^2]$ with the following properties.
\renewcommand{\labelenumi}{(\roman{enumi})}
\begin{enumerate}
\item There is a  unique irreducible component  $C$ on which $f$ has positive degree.  This component meets the diagonal transversally and only at orbifold points.  It is called the {\em handle}.
\item The connected components of the complement of the handle are called the {\em teeth}.  The nodes joining the teeth to the handle are all orbifold points.
\item All ordinary marked points of $C$ lie on the handle.
\end{enumerate}
\end{definition}
This is similar in appearance to Kollar's~\cite{Kol}, but we have additional conditions concerning the marked points and transversality to the diagonal.  

The comb curves form a locally closed substack of $\overline{M}_n(d,g)$ which we denote $U_n(d, g)$.  In fact, as we will see below, they form an {\em open} substack.

If $C$ is a comb curve whose orbifold points are labelled by the set $[2g + 2] = \set{1, 2, \ldots, 2g + 2}$, then we may associate to it a partition of $[2g + 2]$ according to how the orbifold points are distributed among the teeth.  If an orbifold point $x_i$ is on the handle, the corresponding partition includes the singleton set $\set{i}$.  In any such partition, the size of each part must be odd, since the nodes joining the teeth to the handle are orbifold points and there must be an even number of orbifold points on any irreducible component of $C$.

For each partition ${\bf h}$ of $[2g + 2]$ into sets of odd order, let $U_n(d, {\bf h})$ be the moduli space of comb curves with that partition type.   Write $2 g({\bf h}) + 2$ for the number of parts in the partition and label them $h_1, \ldots, h_{2 g({\bf h}) + 2}$.  Then, clearly,
\begin{equation*}
  U_n(d, {\bf h}) = M_n(d, g({\bf h})) \times \prod_{\substack{1 \leq i \leq 2 g({\bf h}) + 2 \\ \# h_i \not= 1}} \overline{M}_0\left(B S_2, \frac{h_i - 1}{2}\right) .
\end{equation*}
For each partition ${\bf h}$, there is a natural map $r : U_n(d,{\bf h}) \rightarrow M_n^{\circ}(d, g({\bf h}))$ sending a comb curve to its handle.  We shall prove
\begin{theorem} \label{thm:mod:1}
  Each $U_n(d, {\bf h})$ is smooth and non-empty and the embedding $U_n(d, {\bf h}) \rightarrow \overline{M}_n(d, g)$ is open.  The virtual degree of $r : U_n(d, {\bf h}) \rightarrow M^\circ_n(d,g({\bf h}))$ is $(-\frac{1}{4})^{g - g({\bf h})}$.
\end{theorem}

By the {\em virtual degree}, we mean the degree of the push-forward of the virtual fundamental class on a fiber, so the theorem asserts
\begin{equation*}
        r_\ast [ U_n(d, {\bf h}) ]^{\rm vir} = \left( -\frac{1}{4} \right)^{g - g({\bf h})} [ M^\circ_n(d, g({\bf h})) ] .
\end{equation*}

The proof occupies the next 3 sections.  In Sections~\ref{proper-intersection} and~\ref{comps-in-diagonal}, we will prove that each embedding $U_n(d, {\bf h}) \rightarrow \overline{M}_n(d,g)$ is open.  In Section~\ref{proof-completion} we reduce the virtual degree calculation to a Hurwitz--Hodge integral that is computed in~\hodgeintegralreference.

\modulisection{Proper intersection with the diagonal}
\label{proper-intersection}

Suppose $f : C \rightarrow [\Sym^2 \PP^2]$ is a representable morphism and $C$ meets $\Delta$ properly: no component of $C$ is carried into $\Delta$ by $f$.  Let $C'$ be a first-order deformation of $C$ with ideal ${\cal O}_C$ and let $P'$ be a closed substack of $C'$.   Assume there are morphisms $f : C \rightarrow [\Sym^2 \PP^2]$ and $g : P' \rightarrow [\Sym^2 \PP^2]$ agreeing on $P = C \fprod[C'] P'$ and  that $P$ is regularly embedded in $C$.  We have a commutative diagram of solid arrows,

\begin{equation*} \label{lifting-problem} \xymatrix{
	P \ar[r] \ar[d] & P' \ar[d] \ar@/^14pt/[ddr]^g \\
	C \ar[r] \ar@/_14pt/[rrd]_f & C' \ar@{-->}[dr]^{f'} \\
	& & [\Sym^2 \PP^2] ,
} \end{equation*}
and we search for a dashed arrow, $f'$ rendering the diagram commutative.  There is an obstruction to the existence of $f'$ in the cohomology group
\begin{equation*}
  H^1(C, f^\ast T [\Sym^2 \PP^2] \tensor {\cal O}_C(-P))
\end{equation*}
where ${\cal O}_C(-P)$ is the ideal sheaf of $P$ in $C$.  If this obstruction is zero then the lifts are a principal homogeneous space under 
\begin{equation*}
  H^0(C, f^\ast T [\Sym^2 \PP^2] \tensor {\cal O}_C(-P)).
\end{equation*}

We are interested in this problem in particular when $P$ is the preimage of a regularly embedded point of the coarse moduli space.  In this case, we can calculate the expected dimension.

\begin{proposition} \label{chi-of-T-of-minus-P}
  Let $r : C \rightarrow \overline{C}$ be the coarse moduli space.  Suppose that $Q$ is a regularly embedded point of $\overline{C}$ and $P = r^{-1}(Q)$.  Then
  \begin{equation*}
    \chi(f^\ast T [\Sym^2 \PP^2] \tensor {\cal O}(-P)) = 3d - 2g - 2
  \end{equation*}
\end{proposition}
\begin{proof}
  We have already done most of the work in \ref{sect:vdim-of-sym2p2}.  The only difference is to compute the degree of $T [\Sym^2 \PP^2] \tensor {\cal O}(-P)$, which is
  \begin{equation*}
    \int_C T [\Sym^2 \PP^2] \tensor {\cal O}(-P) = \int_C T[\Sym^2 \PP^2] - \rank(T [\Sym^2 \PP^2]) \length(P) = 3d - 4
  \end{equation*}
  since $C \rightarrow \overline{C}$ has degree $1$.  Note that the age of $ T [\Sym^2 \PP^2] \tensor {\cal O}(-P)$ is the same as the age of $T [\Sym^2 \PP^2]$ on a fiber because ${\cal O}(-P)$ is pulled back from the coarse moduli space.  The age was determined in Section~\ref{sect:vdim-of-sym2p2} to be $2g + 2$, so orbifold Riemann--Roch gives
  \begin{equation*} \begin{split}
    \chi(f^\ast T [\Sym^2 \PP^2] \tensor {\cal O}(-P))&  = 4 + \deg(f^\ast T [\Sym^2 \PP^2] \tensor {\cal O}(-P)) - \age(f^\ast T [\Sym^2 \PP^2] \tensor {\cal O}(-P)) \\
    &  = 3d - 2g - 2 .
    \end{split}  \end{equation*}
\end{proof}

\begin{proposition} \label{h1-of-T}
  Let $f : C \rightarrow [\Sym^2 \PP^2]$ be a representable morphism from an orbifold pre-stable curve to $[\Sym^2 \PP^2]$.  Assume that $C$ meets the diagonal properly.  Then $H^1(C, f^\ast T [\Sym^2 \PP^2])$ vanishes.  If $P \subset C$ is a closed substack then 
\begin{equation*}
  H^0(C, f^\ast T [\Sym^2 \PP^2]) \tensor H^0(P, {\cal O}_P) \rightarrow H^0(P, f^\ast T [\Sym^2 \PP^2] \rest{P})
\end{equation*}
is surjective.
\end{proposition}
\begin{proof}
  This does not follow immediately from \ref{glob-gen-of-tan} because $f \rest{P}$ need not be a closed embedding.  However, we do know that $f$ is generically a closed embedding because $f$ is representable and it does not carry any component of $C$ into the diagonal.

  Write $T = f^\ast T [\Sym^2 \PP^2]$ and let $V = H^0(C, T) \tensor {\cal O}_C$.  Consider the sequence
  \begin{equation*}
    0 \rightarrow K \rightarrow V \rightarrow T \rightarrow M \rightarrow 0
  \end{equation*}
  where $K$ and $M$ are the kernel and cokernel of $V \rightarrow T$, respectively.  Since $V \rightarrow T$ is surjective on $U$, it follows that $M \rest{U} = 0$, and since $C$ is a curve, $M$ is therefore supported in dimension $0$.

  The spectral sequence computing the cohomology of the above sequence must converge to zero because the sequence is exact.  But the $E_1$ term is
  \begin{equation*}\xymatrix{
      0 \ar[r] & H^1(K) \ar[r] & 0 \ar[r] & H^1(T) \ar[r] & 0 \ar[r] & 0 \\
      0 \ar[r] & H^0(K) \ar[r] & H^0(V) \ar[r]^{\cong} & H^0(T) \ar[r] & H^0(M) \ar[r] & 0 .
    }
  \end{equation*}
  (We have used $H^1(V) = 0$ because $C$ has genus $0$ and $H^1(M) = 0$ because $M$ is supported in dimension $0$.)  The sequence degenerates at the $E_2$ term, which is
  \begin{equation*} \xymatrix{
      0 \ar[drr] & H^1(K) \ar[drr] & 0 \ar[drr] & H^1(T) & 0 \\
      0 & 0 & 0 & 0 & H^0(M) .
    }
  \end{equation*}
  This implies $H^1(K) = H^1(T) = H^0(M) = 0$.

  Now, let $P \subset C$ be any closed substack.  We have an exact sequence
  \begin{equation*}
    V \rest{P} \rightarrow T \rest{P} \rightarrow M \rest{P} \rightarrow 0 
  \end{equation*}
  and a surjection $M \rightarrow M \rest{P}$.  But $M$ and $M \rest{P}$ are supported in dimension zero.  Since we are working in characteristic zero, taking global sections on a DM stack of dimension zero is exact, so $H^0(P, M \rest{P})$ is a quotient of $H^0(C, M)$, and we have just seen that $H^0(C, M) = 0$.  Therefore $H^0(P, T \rest{P})$ is a quotient of 
  \begin{equation*}
    H^0(P, V \rest{P}) = H^0(P, H^0(C, T) \tensor {\cal O}_P) = H^0(C, T) \tensor H^0(P, {\cal O}_P) .
  \end{equation*}
\end{proof}

\begin{corollary} \label{glob-gen-of-T}
  Continue to assume that $f : C \rightarrow [\Sym^2 \PP^2]$ is a representable morphism meeting the diagonal properly.  Let $\pi : C \rightarrow \overline{C}$ be the coarse moduli space.  Whenever $P = \pi^{-1}(Q)$ for some closed point $Q \in \overline{C}$, the map
  \begin{equation*}
    H^0(C, f^\ast T [\Sym^2 \PP^2]) \rightarrow H^0(P, f^\ast T [\Sym^2 \PP^2] \rest{P})
  \end{equation*}
  is surjective and $H^1(C, f^\ast T [\Sym^2 \PP^2] (-P))$ vanishes.
\end{corollary}
\begin{proof}
  Write $T = f^\ast T [\Sym^2 \PP^2]$ as before.  We have seen that
  \begin{equation*}
    H^0(C, T) \tensor H^0(P, {\cal O}_P) \rightarrow H^0(P, T \rest{P})
  \end{equation*}
  is surjective in \ref{h1-of-T}.  The first claim follows from the isomorphism $H^0(P, {\cal O}_P) \cong \CC$.

  For the second claim,   the long exact sequence of cohomology associated to the short exact sequence,
  \begin{equation*}
    0 \rightarrow T(-P) \rightarrow T \rightarrow T \rest{P} \rightarrow 0,
  \end{equation*}
  yields the exact sequence,
  \begin{equation*}
    H^0(C, T) \rightarrow H^0(C, T \rest{P}) \rightarrow H^1(C, T(-P)) \rightarrow H^1(C, T) .
  \end{equation*}
  We have just proven that the first arrow is surjective and we saw in Proposition~\ref{h1-of-T} that $H^1(C, T) = 0$, whence $H^1(C, T(-P)) = 0$.
\end{proof}
  
We have already seen that the obstruction to finding a solution to the lifting problem \eqref{lifting-problem} is a class in $H^1(C, T(-P))$.  If $f : C \rightarrow [\Sym^2 \PP^2]$ meets the diagonal properly then this is the zero vector space: in this case, every such problem has a solution.  Moreover, the space of solutions is a torsor under $H^0(C, T(-P))$.  It has the expected dimension, $3d - 2g - 2$, where $d = \deg(f)$ and $2g + 2$ is the number of orbifold points on $C$.

\begin{corollary} \label{eval-is-smooth}
  Let ${\frk M}_n' \subset {\frk M}_n([\Sym^2 \PP^2])$, $n = 0,1$ be the moduli space parameterizing orbifold pre-stable maps $f : C \rightarrow [\Sym^2 \PP^2]$ with $n$ ordinary marked points and an arbitrary number of orbifold points, and such that no irreducible component of $C$ is carried by $f$ into the diagonal.  The maps
  \begin{equation*} \begin{split} 
      {\frk M}_0' & \rightarrow {\frk M}_0(B \mu_2) \\
      {\frk M}_1' & \rightarrow {\frk M}_1(B \mu_2) \times [\Sym^2 \PP^2] \\
    \end{split}
  \end{equation*}
  are smooth of the expected relative dimensions $3d - 2g + 2$ and $3d - 2g - 2$, respectively.  In particular, the ${\frk M}_n'$, $n = 0,1$ are smooth.
\end{corollary}
\begin{proof}
  In the case $n = 1$, the fiber of the relative obstruction bundle at $(C, f, P)$ is the vector space $H^1(C, f^\ast T [\Sym^2 \PP^2](-P))$, which we have just seen is zero.  If $n = 0$, the relative obstruction bundle is $H^1(C, f^\ast T [\Sym^2 \PP^2])$ which we have also seen to be zero.  The smoothness of the spaces ${\frk M}_n$, then follows from the smoothness of ${\frk M}_n(B \mu_2)$ and of $[\Sym^2 \PP^2]$.  The relative dimensions were determined in Section~\ref{sect:vdim-of-sym2p2} and Proposition~\ref{chi-of-T-of-minus-P}.
\end{proof}

Our next goal is to show that a generic curve meeting the diagonal properly in fact meets it transversally and only at orbifold points.

\begin{proposition} \label{intersection-with-delta}
  Let $f : C \rightarrow [\Sym^2 \PP^2]$ be a representable morphism meeting the diagonal properly.  Let $P$ be the pre-image of a point of the coarse moduli space of $C$ such that $f \rest{P}$ factors through the diagonal.  Fix a small extension $C'$ of $C$.  Then there exists an extension $f'$ of $f$ to $C'$ rendering the solid arrows the the diagram,
  \begin{equation*} \xymatrix{
      P \ar[r] \ar@{-->}[d] & C \ar[r]^f \ar[d] & [\Sym^2 \PP^2] \\
      P' \ar@{-->}[r] & C' \ar[ur]_{f'} ,
    }
  \end{equation*}
  commutative, but such that there is no extension $P'$ of $P$ in $C'$ (the dashed arrows) with $f' \rest{P'}$ factoring through the diagonal.
\end{proposition}
\begin{proof}
  This is a simple dimension argument.  We have seen that for $f$ and $C$ fixed, the space of solid diagrams as above is a principal homogeneous space under $H^0(C, f^\ast T [\Sym^2 \PP^2])$, which has dimension $3d - 2g + 2$.

  On the other hand, consider the space dashed arrows completing the diagram
  \begin{equation*} \xymatrix{
      P \ar@{-->}[r]^a \ar[d] & P' \ar@{-->}[r]^c \ar@{-->}[d]_b & \Delta \ar[d] \\
      C \ar[r] & C' \ar@{-->}[r]^<>(0.5)d & [\Sym^2 \PP^2] 
    }
  \end{equation*}
  while the solid arrows remain fixed.  The space of choices of $a$ and $b$ is a principal homogeneous space under $H^0(P, T_P C)$, which is $1$-dimensional.  Once $a$ and $b$ are fixed, the space of choices for $c$ is a torsor under $H^0(P, f^\ast T \Delta)$, hence of dimension $2$.  Finally, with $a$, $b$, and $c$ all fixed, the choices for $d$ are a torsor under $H^0(C, f^\ast T [\Sym^2 \PP^2](-P))$, giving $3d - 2g - 2$ dimensions of freedom.  Adding these together, we get $3d - 2g + 1$, which is smaller than the dimension of the space of choices for $f'$ without constraining a small extension of $P$ to lie in $\Delta$.  Therefore there are extensions $f'$ in which no small extension factors through $\Delta$.
\end{proof}

\begin{corollary} \label{transversality}
  If $f : C \rightarrow [\Sym^2 \PP^2]$ is generic and meets the diagonal properly then it meets the diagonal transversally and only at orbifold points.
\end{corollary}
\begin{proof}
  A fixed $f : C \rightarrow [\Sym^2 \PP^2]$ meeting the diagonal properly has only finitely many intersection points with the diagonal.  Since the desired property for any fixed point is an open condition, it will be sufficient to deform any bad points away from the diagonal, one at a time.  The choices of $(C, f)$ vary in a smooth space, so it is sufficient to produce a first-order deformation for a single point, which is accomplished by Proposition~\ref{intersection-with-delta}.
\end{proof}

\modulisection{Components in the diagonal}
\label{comps-in-diagonal}

We have determined that ${\frk M}^{\circ}_n([\Sym^2 \PP^2], d,g)$ is smooth of the expected dimension, $3d + 1$.  Any excess dimension in ${\frk M}_n([\Sym^2 \PP^2], d, g)$ must therefore come from curves with components in the diagonal.  We study those curves now.

We will say that an extension $f' : C' \rightarrow [\Sym^2 \PP^2]$ of $f : C \rightarrow [\Sym^2 \PP^2]$ moves a component $C_1 \subset C$ out of the diagonal if $f \rest{C_1}$ factors through $\Delta$ but the restriction of $f'$ to the first order neighborhood $C'_1$ of $C_1$ in $C'$ does not factor through $\Delta$.  If $f'$ moves $C'$ out of $\Delta$, then the induced map on the normal bundles, $N_{C_1 / C'} \rightarrow f^\ast N_{\Delta / [\Sym^2 \PP^2]}$ must be nonzero.

\begin{proposition} \label{deformation-of-nodes}
  Let $C_0$ be an irreducible component of $C$ with $2g + 2$ orbifold points.  Suppose that $f$ has degree $d$ on $C_0$ and $f \rest{C_0}$ factors through $\Delta$.  Let $(C', f')$ be a first-order deformation of $(C, f)$ that moves $C_0$ out of $\Delta$ and smoothes $k$ of the {\em external nodes} of $C_0$ (the nodes joining $C_0$ to the rest of $C$).  Then $d + k \geq g + 1$.  
\end{proposition}
\begin{proof}
  The map $f'$ induces a homomorphism $N_{C_0/C'} \rightarrow f^\ast N_{\Delta / [\Sym^2 \PP^2]}$ of sheaves on $C_0$.  Since $f'$ moves $C_0$ out of $\Delta$, this homomorphism is nonzero.  We have $N_{C_0 / C'} = {\cal O}_{C_0}(-\sum_{i = 1}^k P_i)$, with the sum taken over nodes smoothed to first order in $C'$.  Since $f_0 = f \rest{C_0}$ factors through $\Delta$, we know that $f^\ast N_{\Delta / [\Sym^2 \PP^2]} \cong \rho_1 \tensor p^\ast \overline{f}_0^\ast T \PP^2$, where $p : C \rightarrow \overline{C}$ is the coarse moduli space and $\overline{f}_0 : \overline{C}_0 \rightarrow \PP^2$ is the map induced by $f$.  We therefore obtain a nonzero section of
  \begin{equation*}
    F = {\cal O}(\sum_{i = 1}^k P_i) \tensor \rho_1 \tensor p^\ast \overline{f}_0^\ast T \PP^2 
  \end{equation*}
  over $C_0$.  Noting that sections of $F$ over $C_0$ are in bijection with sections of $p_\ast F$ over $\overline{C}_0$ we obtain a section of
  \begin{equation*}
    p_\ast F = {\cal O}\left(\frac{k}{2} - \frac{2g + 2 - k}{2} \right) \tensor \overline{f}_0^\ast T \PP^2 = {\cal O}(k - g - 1) \tensor \overline{f}_0^\ast T \PP^2 
  \end{equation*}
  over $\overline{C}_0$.  

  But the Euler sequence implies $\overline{f}_0^\ast T\PP^2$ is a quotient of ${\cal O}(\frac{d}{2})^{\oplus 3}$.  As $\overline{C}_0$ has genus zero, $\overline{f}_0^\ast T \PP^2$ must split into ${\cal O}(a) \oplus {\cal O}(b)$ with $\frac{d}{2} \leq a \leq b \leq d$.  Therefore $\overline{f}_0^\ast T \PP^2 \cong {\cal O}(a + k - g - 1) \oplus {\cal O}(b + k - g - 1)$.  We have argued that this bundle must have a nonzero section.  Since $b \leq d$, this implies that $d + k - g - 1 \geq 0$.
\end{proof}

\begin{corollary} \label{comb-nodes}
  If $f : C \rightarrow [\Sym^2 \PP^2]$ is a comb curve and $C \rightarrow C'$ is a small extension smoothing one of the nodes that joins a tooth of $C$ to the handle, then there is no extension of $f$ to $f' : C' \rightarrow [\Sym^2 \PP^2]$.
\end{corollary}
\begin{proof}
  Let $C_0$ be the handle of $C$ and let $C_1$ be a tooth joined at an orbifold point $P \in C_0$.  Assume for the moment that $C_1$ is irreducible.  The fiber of $N_{C_1 / C'}$ at $P$ is generated by $T_P C_0$.  Since $C_0$ meets the diagonal transversally at $P$, the map $T_P C_0 \rightarrow f^\ast N_{\Delta / [\Sym^2 \PP^2]}$ is nonzero.  Thus any extension of $f$ to $C'$ must determine a nonzero map $N_{C_1 / C'} \rightarrow f^\ast N_{\Delta / [\Sym^2 \PP^2]}$: it must move $C_1$ out of the diagonal.  But $f \rest{C_1}$ has degree zero, so by the proposition, any such extension must smooth at least $g + 1$ nodes of $C_1$, with $2g + 2$ being the number of orbifold points on $C_1$.  By stability, there are at least $2$ orbifold points on $C_1$, so $g + 1 \geq 1$, but there is only one node on $C_1$, hence no such smoothing can exist.

  If $C_1$ is reducible, we proceed by induction on the components of $C_1$, since at least one of the nodes of $C_1$ must be smoothed in this case by the proposition.  Repeating the argument on the branch attached at this node (which must have fewer irreducible components) completes the proof.
\end{proof}

\begin{corollary} \label{comb-curves-open}
  For each partition ${\bf h}$, the embedding  $U_n(d,{\bf h}) \rightarrow \overline{M}_n(d,g)$ is open.
\end{corollary}
\begin{proof}
  The last corollary showed that any small extension of a comb curve is a comb curve, or, in other words, that the embedding $U_n(d,{\bf h}) \rightarrow \overline{M}_n(d,g)$ is smooth, and therefore an open embedding.
\end{proof}

\begin{corollary} \label{comb-curves-disj-union}
  The     locus of comb curves in $\overline{M}_n(d,g)$ breaks into a disjoint union
  \begin{equation*}
    U_n(d,g) = \coprod_{{\bf h}} U_n(d, {\bf h})
  \end{equation*}
  over all partitions ${\bf h}$ of $[2g + 2]$ into subsets of odd orders.
\end{corollary}
\begin{proof}
  By the last corollary, each $U_n(d, {\bf h})$ is open in $U_n(d,g)$.
\end{proof}

\modulisection{The obstruction bundle for comb curves}
\label{proof-completion}
\label{obstruction-bundle}

Unfortunately, the results of the previous sections do not give us a complete understanding of the moduli spaces $\overline{M}_n([\Sym^2 \PP^2],d,g)$.  However, Corollary~\ref{comb-nodes} does give an essentially complete description of the locus of comb curves.

In this section, we will need 
\begin{proposition}[Behrend--Fantechi \cite{BF}, Proposition 5.6]
  If $M$ is a smooth Deligne--Mumford stack with an absolute obstruction theory $E = [E^{-1} \rightarrow E^0]$ then $E^{-1}$ is a vector bundle on $M$ and the virtual fundamental class of $M$ is $c_{top}({E^{-1}}^\vee)$.
\end{proposition}

Applying this to $U_n(d,g)$, we see that its virtual fundmental class is the top Chern class of its {\em absolute} obstruction bundle, which we denote $\Obs(C, f)$.  We may gain access to this bundle via the tangent--obstruction sequence, a fragment of which is
\begin{equation*}
  \Def(C) \rightarrow \Obs(f) \rightarrow \Obs(C, f) \rightarrow 0.
\end{equation*}

From Corollary~\ref{comb-nodes} we know that no node of a comb curve attaching a tooth to the handle can be smoothed.  However, it is a straightforward consequence of Corollary~\ref{eval-is-smooth} that $f$ can be extended to any first-order  deformation $C'$ of $C$ that does not smooth any node joining a tooth to the handle.  Thus the image of the map
\begin{equation*}
  \Def(C) \rightarrow \Obs(f)
\end{equation*}
is the vector space parameterizing deformations of the nodes that join the teeth to the handle (that is, deformations of $C$ modulo deformations that do not smooth those nodes to first-order).  This space is $\sum \pi_\ast \left( T_{P_i} C_0 \tensor T_{P_i} C_i \right)$, the sum being taken over the teeth $C_i$, with $P_i$ being the node joining $C_i$ to $C_0$, and $\pi$ being the projection from $C$ to $U_n(d,g)$.

Combining this with the tangent--obstruction sequence gives a short exact sequence,
\begin{equation} \label{obs-seq}
  0 \rightarrow \sum \pi_\ast \left( T_{P_i} C_0 \tensor T_{P_i} C_i \right) \rightarrow \Obs(f) \rightarrow \Obs(C, f) \rightarrow 0.
\end{equation}
The middle term can be computed explicitly.
\begin{lemma}
  If $f : C \rightarrow [\Sym^2 \PP^2]$ is a comb curve with teeth $C_i$, $i = 1, \ldots, k$ then 
  \begin{equation*}
    H^1(C, f^\ast T [\Sym^2 \PP^2]) = \sum_{i = 1}^k H^1(C_i, f^\ast T [\Sym^2 \PP^2] \rest{C_i}) .
  \end{equation*}
\end{lemma}
\begin{proof}
  Let $\nu : C^\nu = \coprod_{i = 0}^k C_i \rightarrow C$ be the normalization of the nodes $P_i$, $i = 1, \ldots, k$ that join the teeth to the handle (taking $C_0$ to be the handle).  Let $T = f^\ast T [\Sym^2 \PP^2]$.  Then the normalization sequence on $C$ yields the exact sequenece,
  \begin{equation*}
    \sum_{i = 0}^k H^0(C_i, T \rest{C_i}) \rightarrow \sum_{i = 1}^k H^0(P_i, T \rest{P_i}) \rightarrow H^1(C, T) \rightarrow \sum_{i = 0}^k H^1(C_i, T \rest{C_i}) \rightarrow 0 .
  \end{equation*}
  The first arrow is surjective by Corollary~\ref{glob-gen-of-T}, so 
  \begin{equation*}
        H^1(C, T) = \sum_{i = 0}^k H^1(C_i, T \rest{C_i}) .
  \end{equation*}     
  But $H^1(C_0, T \rest{C_0}) = 0$ by Proposition~\ref{h1-of-T}, since the handle meets the diagonal properly, whence the lemma.
\end{proof}

Now we have
\begin{equation*}
  \Obs(f) = H^1(C, f^\ast T [\Sym^2 \PP^2]) = \sum_{i = 1}^k H^1(C_i, f^\ast T [\Sym^2 \PP^2] \rest{C_i}) .
\end{equation*}
But $C_i \rightarrow [\Sym^2 \PP^2]$ factors through a point in the diagonal.  Hence $f^\ast T [\Sym^2 \PP^2] \rest{C_i} \cong \rho_0^{\oplus 2} \oplus \rho_1^{\oplus 2}$ where $\rho_0$ and $\rho_1$ are the trivial and non-trivial representations of $\mu_2$, respectively.  Thus,
\begin{equation*}
  \Obs(f) = \sum_{i = 1}^k H^1(C_i, \rho_0^{\oplus 2} \oplus \rho_1^{\oplus 2}) .
\end{equation*}
But $\rho_0 \cong {\cal O}_{C_i}$ has no higher cohomology because $C_i$ has genus zero.  We are left with
\begin{equation*}
  \Obs(f) = \sum_{i = 1}^k H^1(C_i, \rho_1)^{\oplus 2} .
\end{equation*}
Defining ${\bf E}_i^\vee = R^1 \pi_\ast \rho_1$ for $\pi : C \rightarrow U_n(d,g)$ the universal curve, we have therefore proven that there is an exact sequence
\begin{equation*}
  0 \rightarrow \sum_{i = 1}^k \pi_\ast \left( T_{P_i} C_i \tensor T_{P_i} C_0 \right) \rightarrow \sum_{i = 1}^k {\bf E}_i^\vee \oplus {\bf E}_i^\vee \rightarrow \Obs(C, f) \rightarrow 0  .
\end{equation*}
To keep the notation readable, let us now write $T_i = T_{P_i} C_i$ for $i \not= 0$ and $T'_i = T_{P_i} C_0$.  The above sequence determines the total Chern class of $\Obs(C,f)$ on $U_n(d,g)$ to be
\begin{equation*}
  \prod_{i = 1}^{2 g({\bf h}) + 2} \frac{c({\bf E}_i^\vee)^2}{c(\pi_\ast(T_i \tensor T'_i))} .
\end{equation*}
In order to complete the proof of Theorem~\ref{thm:mod:1} we must integrate this Chern class on a fiber of the map $r : U_n(d, {\bf h}) \rightarrow M^\circ_n(d,g({\bf h}))$.

Let $h_1, \ldots, h_{2g({\bf h}) + 2}$ be the sets in the partition ${\bf h}$ and let $2 g_i + 1$ be the number of elements in $h_i$.  On a fiber $r$, $T'_i$ is isomorphic to $\rho_1$, so the integral becomes
\begin{equation*}
  \int_{\prod_i \overline{M}_0(B S_2, g_i)} \prod_i \frac{c({\bf E}_i^\vee)^2}{c(\pi_\ast (T_i \tensor \rho_1))} = \prod_i \int_{\overline{M}_0(B S_2, g_i)} \frac{c({\bf E}_i^\vee)^2}{1 - \psi_1} .
\end{equation*}
The integral under the product is evaluated in \hodgeintegralreference.  Its value is $(-\frac{1}{4})^{g_i}$.  Taking the product over $i$ and noting that $g({\bf h}) + \sum g_i = g$, we obtain $(- \frac{1}{4})^{g - g({\bf h})}$ for the virtual degree.  This completes the proof of Theorem~\ref{thm:mod:1}.

\subsection{The evaluation map}
\label{evaluation}

Let $M = \overline{M}_n(d, g)$.  There is an evaluation map $M \rightarrow [\Sym^2 \PP^2]^n$.  The goal of this section is to estimate the dimension of the image of this evaluation map.

\begin{lemma}
  The image of the evaluation map $\overline{M}_n([\Sym^2 \PP^2], d, g) \rightarrow [\Sym^2 \PP^2]^n$ has dimension at most $3d + 1 + n$.  The image of the locus parameterizing curves that have more than one component with positive degree is of strictly smaller dimension.
\end{lemma}
\begin{proof}
  The proof is by induction on the number of components with positive degree and the number of components with image in the diagonal.  To be slightly more precise, these properties define a finite stratification of $\overline{M}_n([\Sym^2 \PP^2], d, g)$ and we prove the result for one stratum at a time.

  Consider first the open stratum $M_n([\Sym^2 \PP^2], d, g)$ parameterizing maps from irreducible curves to $[\Sym^2 \PP^2]$.  This is the union of the closed substack $T_1(n,d,g) = \overline{M}_n(\Delta, d, g)$ and its complement $T_2(n,d,g) = \overline{M}_n([\Sym^2 \PP^2], d, g) - T_1(n,d,g)$.  Now, the evaluation map $T_1(n,d,g) \rightarrow [\Sym^2 \PP^2]^n$ factors through $\Delta^n$.  Furthermore, the composition of the evaluation map with $\Delta^n \rightarrow (\PP^2)^n$ factors through $M_{0,n}(\PP^2, \frac{d}{2})$, as in the diagram
\begin{equation*} \xymatrix{
    T_1(n,d,g) \ar[r] \ar[d] & \Delta^n \ar[r] \ar[d] & [\Sym^2 \PP^2]^n \\
    \overline{M}_{0,n}(\PP^2, \frac{d}{2}) \ar[r] & \left( \PP^2 \right)^n .
} \end{equation*}
Now, $\Delta^n \rightarrow (\PP^2)^n$ is a gerbe, so the dimesnion of the image of $T_1(n,d,g)$ in $\Delta^n$ coincides with the dimension of the image of $\overline{M}_{0,n}(\PP^2, \frac{d}{2})$ in $\left( \PP^2 \right)^n$.  This latter number is bounded by $\dim \overline{M}_{0,n}(\PP^2, \frac{d}{2}) = \frac{3d}{2} - 1 + n$.  This proves the lemma for the stratum $T_1(n,d,g)$.

For $T_2(n,d,g)$, we may refer to \ref{vdim-calcs}, which implies $\dim T_2(n,d,g) = 3d + 1 + n$.  Thus the lemma also holds for $T_2(n,d,g)$.  

  Assume now that the conclusion of the lemma holds for the open substacks 
  \begin{equation*}
    V^{(k)}_n(d,g) \subset \overline{M}_n([\Sym^2 \PP^2], d, g)
  \end{equation*}
  parameterizing orbifold stable maps to $[\Sym^2 \PP^2]$ with at most $k$ irreducible components.  Then we may obtain $V^{(k + 1)}_{n''}(d'', g'')$ as the union of the stacks
  \renewcommand{\labelenumi}{\Roman{enumi}.}
  \begin{enumerate}
  \item $\displaystyle      V^{(k)}_{n + 1}(d,g) \fprod[{[\Sym^2 \PP^2]}] T_1(n' + 1, d',g')$
  \item $\displaystyle       V^{(k)}_{n + 1}(d,g) \fprod[{[\Sym^2 \PP^2]}] T_2(n' + 1, d',g')$
  \item $\displaystyle       V^{(k)}_n(d,g) \fprod[\PP^2] T_1(n', d',g')$
  \item $\displaystyle      V^{(k)}_n(d,g) \fprod[\PP^2] T_2(n', d',g')$,
  \end{enumerate}
  the union being taken over all partitions $n'' = n + n'$, $d'' = d + d'$, and $g'' = g + g'$.

  These parameterize, respectively, 
  \begin{enumerate}
  \item curves with at most $k$ components joined at an ordinary point to a smooth curve in $\Delta$, 
  \item curves with at most $k$ components joined at an ordinary point to a smooth curve meeting $\Delta$ properly, 
  \item curves with at most $k$ components joined at an orbifold point to a smooth curve in $\Delta$, and
  \item curves with at most $k$ components joined at an orbifold point to a smooth curve meeting $\Delta$ properly.  
  \end{enumerate}
  Note that we only need to consider adjoining a new component at a single point because we are only working with genus zero curves.

  Before considering these cases individually, we note that by attaching a component of degree $d' = 0$, we can increase the dimension of the image by at most $1$: when $n' \geq 0$ the dimension increases by $1$; otherwise it does not increase at all.  Therefore we may assume $d' > 0$ below.

  \renewcommand{\labelenumi}{\sc Case \Roman{enumi}.}
  \begin{enumerate}
  \item The evaluation map on $V_{n + 1}(d, g) \fprod[{[\Sym^2 \PP^2]}] T_1(n', d',g')$ factors through $[\Sym^2 \PP^2]^n \times \Delta^{n'}$.  Composing with the map to the coarse moduli space, $[\Sym^2 \PP^2]^n \times \Delta^{n'} \rightarrow (\Sym^2 \PP^2)^n \times (\PP^2)^{n'}$ does not change the the dimension of the image but the composed evaluation map factors through
    \begin{equation*}
      \image \left( V_{n + 1}(d, g) \rightarrow [\Sym^2 \PP^2]^{n+1} \right) \fprod[\Sym^2 \PP^2] \overline{M}_{0, n' + 1}\left(\PP^2, \frac{d'}{2}\right) .
    \end{equation*}
    (Forgetting orbifold marked points does not destabilize the curve because we have assumed it has positive degree.)

    Now, the evaluation map at the attaching point, $\overline{M}_{0, n'+1}(\PP^2, d') \rightarrow \PP^2$ , is smooth by \ref{eval-is-smooth} so the fiber product above has the expected dimension.  By the inductive hypothesis, it is at most
    \begin{equation*}
      (3d + 1 + n + 1) + \left(\frac{3d'}{2} - 1 + n' + 1\right) - 2 = 3(d + d') + n + n' - \frac{3d'}{2} < 3d'' + 1 + n''
    \end{equation*}
    which completes the induction in this case.

    \item Now consider the evaluation map on $V_{n + 1}(d,g) \fprod[{[\Sym^2 \PP^2]}] T_2(n' + 1, d',g')$.  In this case, the evaluation map factors through
      \begin{equation*}
        \image\left(V_{n + 1}(d, g) \rightarrow [\Sym^2 \PP^2]^{n + 1}\right) \fprod[{[\Sym^2 \PP^2]}] T_2(n' + 1, d', g')
      \end{equation*}
      The map, $T_2(n' + 1,d',g') \rightarrow [\Sym^2 \PP^2]$, that evaluates at the attaching point is smooth be \ref{eval-is-smooth} so the fiber product above has the expected dimension.  By the inductive assumption, it is bounded by
      \begin{equation*}
        (3d + 1 + n + 1) + (3d' + 1 + n' + 1) - 4 = 3(d + d') + (n + n') < 3d'' + 1 + n''.
      \end{equation*}
      This completes the induction in this case.

    \item This is almost exactly the same as {\sc Case I}.
      
    \item The evaluation map on $V_n(d, g) \fprod[\PP^2] T_2(n', d',g')$ factors through
        \begin{equation*}
          \image\left(\overline{M}_n([\Sym^2 \PP^2], d, g) \rightarrow \left( \PP^2 \right)^n \right) \fprod[\PP^2] T_2(n', d', g')
        \end{equation*}
        and the evaluation map at an orbifold point $T_2(n', d', g') \rightarrow \PP^2$ is smooth, so the fiber product above has the expected dimension which is bounded by
        \begin{equation*}
          (3 d + 1 + n) + (3 d' + 1 + n') - 2 = 3 (d + d') + (n + n') < 3d'' + 1 + n''.
        \end{equation*}
        This completes the induction, and the proof.
  \end{enumerate}
\end{proof}

Let $V$ be the locus of curves in $\overline{M}_n([\Sym^2 \PP^2], d, g)$ having a single component with positive degree that does not map into $\Delta$.  By the proposition, the image of the complement of $V$ in $\overline{M}_n(d,g)$ in $[\Sym^2 \PP^2]^n$ has dimension strictly smaller than $3d + 1 + n$.  This will permit us to restrict attention to $V$ for our enumerative applications.  However, it will be advantageous to restrict attention still further to the locus comb curves $U_n(d,g) \subset V$ using 

\begin{proposition}\label{dim-ev-image}
  The image of $U_n(d,g)$ under the evaluation map
  \begin{equation*}
    \overline{M}_n(d,g) \rightarrow [\Sym^2 \PP^2]^n
  \end{equation*}
  has dimension $3d + 1 + n$ and the image of the complement of $U_n(d,g)$ has strictly smaller dimension.
\end{proposition}
\begin{proof}
In view of the lemma, it is sufficient to show that the image of $V - U_n(d,g)$ has dimension strictly smaller than $3d + 1 + n$.  There are four reasons $(C, f) \in V$ may fail to be in $U_n(d,g)$: either a marked point appears on a tooth, a tooth is joined to the handle at an ordinary point, or the handle meets the diagonal at an ordinary point, or the handle meets the diagonal non-transversally at an orbifold point.  

  First, consider the map $V \rightarrow M = M_{n'}^\circ([\Sym^2 \PP^2], d, g')$ which sends $C$ to its handle.  We know by \ref{eval-is-smooth} that $M$ is smooth, and a generic point of $M$ corresponds to a curve in $[\Sym^2 \PP^2]$ that meets $\Delta$ transversally and only at orbifold points.  Thus the locus of curves in $M$ that fail to have these properties has dimension strictly smaller than $\dim M = 3d + 1 + n'$.  

  Note first that $n' < n$ unless all of the marked points of $C$ are on the handle.  The argument of the last paragraph shows that $C$ does not meet the diagonal except at orbifold points and that it must meet the diagonal transversally there.  Finally, since all of the marked points are on the handle, the stability of $(C,f)$ implies that any tooth must contain an orbifold point.  But then the corresponding node maps into the diagonal, hence it is an orbifold point.
\end{proof}

%% file: gw-enum-2008.07.24.tex
s
In this section, we will relate the Gromov--Witten invariants of $[\Sym^2 {\bf P}^2]$ to the enumeration of hyperelliptic curves in ${\bf P}^2$.  In some cases where the enumerative geometry of hyperelliptic curves is simple, this will enable us to compute Gromov--Witten invariants.

\subsection{Notation}

Let $\phi$ be a class in $A_\ast([\Sym^2 \PP^2]^n \times (\PP^2)^{2g + 2})$.  We write
\begin{equation*}
  \GW{\phi}{(d,g)} = \int_{[\overline{M}_n(d,g)]^{\rm vir}} e^\ast(\phi)
\end{equation*}
where $e : \overline{M}_n(d,g) \rightarrow \overline{I} [\Sym^2 \PP^2]^n$ is the evaluation map.  If $\phi \in A_\ast([\Sym^2 \PP^2]^n)$ and 
\begin{equation*}
p : [\Sym^2 \PP^2]^n \times (\PP^2)^{2g + 2} \rightarrow [\Sym^2 \PP^2]^n
\end{equation*}
is the projection, it is also convenient to write $\GW{\phi}{(d,g)}$ instead of $\GW{p^\ast \phi}{(d,g)}$.

We explain the relationship between our notation for Gromov--Witten invariants, $\GW{\hphantom{-}}{(d,g)}$ and the notation $\GW{\hphantom{-}}{d}$ used by Abramovich, Graber, and Vistoli~\cite{AGV}, as the latter will be used in Sections~\ref{gw} and~\ref{crc}.  Let $\set{\phi_i}$ be a collection of homogeneous elements of $A_\ast(\overline{I} [\Sym^2 \PP^2])$ such that each $\phi_i$ comes either from the twisted sector or the untwisted sector.  Let $2g + 2$ be the number of the $\phi_i$ that come from the twisted sector.  Then
\begin{equation*}
\GW{\phi_1, \ldots, \phi_n}{d} = \GW{\phi_1,\ldots,\phi_n}{(d,g)} 
\end{equation*}
for all $d$.  Conversely, we have
\begin{equation*}
\GW{\phi_1, \ldots, \phi_n}{(d,h)} = \begin{cases} \GW{\phi_1, \ldots, \phi_n, \gamma^{\tensor (2h - 2g)}}{d} & h \geq g \\ 0 & h < g . \end{cases}
\end{equation*}
Thus the $\GW{\hphantom{-}}{(d,g)}$ and $\GW{\hphantom{-}}{d}$ package the same information in different ways.

\subsection{The degree $0$ invariants}
\label{deg-0}

We will begin by computing the degree zero invariants of $[\Sym^2 \PP^2]$ using the Chow rings of $[\Sym^2 \PP^2]$ and $\PP^2$ and the calculation by Faber and Pandharipande~\cite{FP} of the hyperelliptic Hodge integral, $\int \lambda_g \lambda_{g - 1}$.

\begin{lemma} \label{deg-0-untwisted-insert}
  Suppose $\phi_1 \in A_\ast(\overline{I} [\Sym^2 \PP^2])$ is a class in the {\em untwisted} sector.  If $n > 3$ or $g > 0$, then
  \begin{equation*}
    \left< \phi_1, \ldots, \phi_n \right>_{(0,g)} = 0 .
  \end{equation*}
\end{lemma}
\begin{proof}   
  By linearity of the Gromov--Witten invariants, we can assume that each $\phi_i \in A_\ast(\overline{\Omega}_j)$ for some $j$.  If $n'$ is the number of $\phi_i$ coming from the untwisted sector, then $n' \geq 1$ and the above invariant is computed on the moduli space $\overline{M}_{n'}([\Sym^2 \PP^2], 0, g)$, which we abbreviate to $T_{n'}$.  Since $n > 3$ or $g > 0$, there is a forgetful map
  \begin{equation*}
    q : T_{n'} \rightarrow T_{n' - 1} .
  \end{equation*}
  These spaces parameterize degree zero maps, so the evaluation map $e : T_{n'} \overline{I}([\Sym^2 \PP^2])^{n'}$ factors through this forgetful map.  Thus
  \begin{equation*}
    e^\ast \phi = q^\ast \overline{e}^\ast \phi
  \end{equation*}
  for a map $\overline{e}$ defined on $T_{n'-1}$.  Moreover $q$ is smooth of the expected dimension, so $q^\ast [ T_{n'-1} ]^{\rm vir} = [T_{n'}]^{\rm vir}$.  Therefore, 
  \begin{equation*}
    \int e^\ast \phi \cap [T_{n'}]^{\rm vir} = \int q^\ast \left(  \overline{e}^\ast \phi \cap [ T_{n'-1}]^{\rm vir} \right) .
  \end{equation*}
  This must be zero because the fibers of $q$ have positive dimension.  
\end{proof}
  
By the lemma, the only potentially nonzero invariants of degree zero are those with $n = 3$ and $g \leq 0$, and those with all insertions in the twisted sector.  We consider the case where $n = 3$ and $g \leq 0$ case first.

Consider the invariant
\begin{equation*}
  \GW{\phi_1, \phi_2, \phi_3}{(0,g)} .
\end{equation*}
If $g = -1$ then all $\phi_i$ come from $A_\ast ([\Sym^2 \PP^2])$ and this is just an integral on $[\Sym^2 \PP^2]$:
\begin{equation*}
  \GW{\phi_1, \phi_2, \phi_3}{(0,-1)} = \int_{[\Sym^2 \PP^2]} \phi_1 \phi_2 \phi_3 .
\end{equation*}

If $g = 0$ then two of the $\phi_i$ come from the twisted sector --- say $\phi_1$ and $\phi_2$.  In this case, we compute the Gromov--Witten invariant on the moduli space parameterizing degree zero orbifold stable maps to $[\Sym^2 \PP^2]$ with $2$ orbifold marked points and one ordinary marked point.  This moduli space is isomorphic to $\Delta$, with the first two evaluation maps to $\overline{\Omega}_1 \cong \PP^2$ being projection on the coarse moduli space and the third being the inclusion in $[\Sym^2 \PP^2]$.  Therefore, we have
\begin{equation*}
  \left< \phi_1, \phi_2, \phi_3 \right>_{(0,0)} = \int_\Delta q^\ast(\phi_1 \phi_2) i^\ast (\phi_3) = \int_{[\Sym^2 \PP^2]} i_\ast q^\ast (\phi_1 \phi_2) . \phi_3
\end{equation*}
where $q : \Delta \rightarrow \PP^2$ is map to the coarse moduli space and $i : \Delta \rightarrow [\Sym^2 \PP^2]$ is the inclusion.  

This completes the calculation of the degree zero invariants involving an insertion in the untwisted sector.  We are left to consider the invariants where all insertions come from the twisted sector.

If there are no ordinary marked points then the expected dimension of the moduli space of degree zero maps to $[\Sym^2 \PP^2]$ is $1$.  Since $A_1(\overline{\Omega}_1) = A_1(\PP^2) = \QQ \gamma_1$, this means that up to linearity, the only remaining degree zero Gromov--Witten invariant of interest is
\begin{equation*}
\GW{\gamma_1}{(0,g)} =   \GW{\gamma_1, \gamma_0^{\tensor (2g + 1)}}{0} = \int_{[\overline{M}_0([\Sym^2 \PP^2], 0, g)]^{\rm vir}} e_1^\ast (\gamma_1)
\end{equation*}
(recall from \ref{group-structure} that $\gamma_1$ is the class of a line in the twisted sector and $\gamma_0$ is the fundamental class of the twisted sector).  We will evaluate this integral on $\overline{M}_0([\Sym^2 \PP^2], 0, g)$, the moduli space of degree $0$ maps to $[\Sym^2 \PP^2]$ with $2g + 2$ orbifold markings and no other markings, which is naturally identified with
\begin{equation*}
  \overline{M}_0(\Delta, 0, g) = \overline{M}_0(B S_2, g) \times \PP^2
\end{equation*}
(the last moduli space parameterizes orbifold stable maps to $B S_2$ with $2g + 2$ orbifold marked points and no other marked points).    Since $M = \overline{M}_0(B S_2, g) \times \PP^2$  is manifestly smooth over $\overline{M}_0(B S_2, g)$ the virtual class equals the top Chern class of the relative obstruction bundle,
\begin{equation*}
  R^1 \pi_\ast f^\ast T [\Sym^2 \PP^2] ,
\end{equation*}
where 
\begin{equation*} \xymatrix{
    C \ar[r]^<>(0.5)f \ar[d]_\pi & [\Sym^2 \PP^2] \\
    M
  } 
\end{equation*}
is the universal curve over $M$.  We note that $f$ factors through $\Delta \cong B S_2 \times \PP^2$ since $C$ has orbifold points and $f$ has degree zero.  Thus the universal map $f$ factors through a map $g : M \rightarrow B S_2 \times \PP^2$.  Note also that 
\begin{equation*}
        T [\Sym^2 \PP^2] \rest{\Delta} \cong (\rho_0 \boxtimes T \PP^2) \oplus (\rho_1 \boxtimes T \PP^2) = (\rho_0 \oplus \rho_1) \boxtimes T \PP^2
\end{equation*}
 where $\rho_0$ and $\rho_1$ is the trivial and non-trivial representations of $S_2$, repsectively, viewed as line bundles on $B S_2$.  Thus,
\begin{equation*}
    R^1 \pi_\ast f^\ast T [\Sym^2 \PP^2] 
    = R^1 \pi_\ast \pi^\ast g^\ast ((\rho_0 \oplus \rho_1) \boxtimes T \PP^2)
    \cong R^1 \pi_\ast (\rho_0 \oplus \rho_1) \boxtimes T \PP^2 .
\end{equation*}
We must calculate $R^1 \pi_\ast (\rho_0 \oplus \rho_1)$.  Put $\tilde{C} = C \fprod[B S_2] (\rm point)$.  Then $\tilde{C}$ is a family of hyperelliptic curves over $M$.  Let $q : \tilde{C} \rightarrow C$ be the projection.  Then $q_\ast {\cal O}_{\tilde{C}} = \rho_0 \oplus \rho_1$.  Since $q$ is affine, this means that
\begin{equation*}
  R^1 \pi_\ast (\rho_0 \oplus \rho_1) = R^1 (\pi_\ast q_\ast) {\cal O}_{\tilde{C}} = {\bf E}^\vee ,
\end{equation*}
the dual of the hyperelliptic Hodge bundle.

We return to the problem of calculating the Gromov--Witten invariant
\begin{equation*}
  \GW{ \gamma_1 }{(0,g)} = \int_{\overline{M}_0(B S_2, g)} c_{\rm top}({\bf E}^\vee \boxtimes T \PP^2) \gamma_1 .
\end{equation*}

Let $a_1, \ldots, a_g$ be the Chern roots of $\bf E$, let $b_1, b_2$ be the Chern roots of $T {\bf P}^2$, and let $\lambda_i = c_i({\bf E})$.  Then 
\begin{equation*}
c_{\rm top}({\bf E}^\vee \boxtimes T {\bf P}^2) = \prod_i (-a_i + b_1)(-a_i + b_2) = \lambda_g^2 - 3 \lambda_g \lambda_{g - 1} h + 3 (\lambda_g \lambda_{g - 2} +  \lambda_{g - 1}^2) h^2 .
\end{equation*}
since $c(T \PP^2) = 1 + 3 h + 3 h^2$, where $h$ is the hyperplane class on $\PP^2$.

Now we compute
\begin{equation*}
\int_{\overline{M}_0(B S_2, g) \times {\bf P}^2} c_{\rm top}({\bf E}^\vee \boxtimes T {\bf P}^2) \gamma_1 = -3 \int_{\overline{M}_0(B S_2, g)} \lambda_g \lambda_{g - 1}.
\end{equation*}
The last integral was computed by Faber and Pandharipande~(\cite{FP}, Corollary to Proposition~3).  (The number indicated below differs from theirs by a factor of $(2g + 2)!$ since $\overline{M}_0(B S_2, g)$ is the moduli space of hyperelliptic curves with an ordering of the $2g + 2$ branch points.)  It is
\begin{equation*}
\int_{\overline{M}_0(B S_2, g)} \lambda_g \lambda_{g - 1} = \frac{(-1)^{g - 1} (2^{2g} - 1) B_{2g}}{2g}
\end{equation*}
where the $B_n$ are the Bernoulli numbers, i.e., $\displaystyle \frac{z}{e^z - 1} = \sum_{n = 0}^\infty B_n \frac{z^n}{n!}$.  We conclude that
\begin{equation*}
\GW{ \gamma_1 }{(0,g)} = \frac{(-1)^g \: (2^{2g} - 1) \: 3 \: B_{2g}}{2g}.
\end{equation*}

This completes the calculation of the degree $0$ invariants.

\subsection{The $2$-point invariants}
\label{2-point}
\label{deg-1-through-2-points}

We calculate some invariants of the form
\begin{equation} \label{2-point-form}
\GW{\phi_1, \phi_2}{(d,g)} .
\end{equation}
Since we will not need all invariants of the form~\eqref{2-point-form}, I have only included calculations of the few we will need, followed by a few comments about the remaining ones in Section~\ref{comments}.  

Considering the virtual dimension of the moduli space $\overline{M}_n(d,g)$, we see that if \eqref{2-point-form} is nonzero, then
\begin{equation*}
\deg(\phi_1) + \deg(\phi_2) = 3d + 3.
\end{equation*}
But $\deg(\phi_i) \leq 4$ for each $i$, so $d \leq 1$.  The case $d = 0$ was already addressed in the last section, so we are left with $d = 1$.

\paragraph{The case $g = -1$}
\label{2-point-deg-1}

Note that $\overline{M}_2(1, -1)$ is isomorphic to $\overline{M}_{0,2}(\PP^2, 1) \times \PP^2$ and there is a commutative diagram
\begin{equation*} \xymatrix{
\overline{M}_{0,2}(\PP^2,1) \ar[r]^<>(0.5){e_i \times \id} \ar[d] & \PP^2 \times \PP^2 \ar[d] \\
\overline{M}_2(1, -1) \ar[r]^{e_i} & [\Sym^2 \PP^2]
} \end{equation*}
for each evaluation map $e_i$.  We therefore have
\begin{equation*} \begin{split}
\left< \alpha^4, \alpha^2 \right>_{(1,-1)} 
& = \int_{\overline{M}_{0,2}(\PP^2, 1) \times \PP^2} (e_1 \times \id)^\ast (h_1^2 + 2 h_1 h_2 + h_2^2) (e_2 \times \id)^\ast (6 h_1^2 h_2^2) \\
& = \int 6 e_1^\ast(h_1^2) e_2^\ast(h_1^2) \tensor h_2^2 = 6 \left< h_1^2, h_1^2 \right>^{\PP^2}_1 = 6 \\
\left< \alpha^4, \beta \right>_{(1,-1)}
& = \int (e_1 \times \id)^\ast (h_1 h_2) (e_2 \times \id)^\ast (6 h_1^2 h_2^2) = 0 \\
\left< \alpha^3, \alpha^3 \right>_{(1,-1)}
& = \int (e_1 \times \id)^\ast (3 h_1^2 h_2 + 3 h_1 h_2^2) (e_2 \times \id)^\ast (3 h_1^2 h_2 + 3 h_1 h_2^2) \\
& = 9 \int e_1^\ast(h_1^2) e_2^\ast(h_1)^2 \tensor h_2^2  = 9 \left< h_1^2, h_1^2 \right>_1^{\PP^2}  = 9 .
\end{split} \end{equation*}

\paragraph{The case $g = 0$}

A point of $\overline{M}_0(1,0)$ determines a map from a genus zero curve with $2$ marked points to $\PP^2$: let $C \rightarrow [\Sym^2 \PP^2]$ be an orbifold stable map; pulling back to $\PP^2 \times \PP^2$ and composing with the projection of $\PP^2 \times \PP^2$ on the first factor gives a map from a curve of genus zero to $\PP^2$; stabilizing this map gives $p : \overline{M}_0(1, 0) \rightarrow \overline{M}_{0,2}(\PP^2, 1)$.  For any point of $M_{0,2}(\PP^2, 1)$ there is a unique involution of the source curve that fixes the $2$ marked points, so $p$ is birational.  We have a commutative diagram,
\begin{equation*} \xymatrix{
\overline{M}_0(1, 0) \ar[r]^{(e_1, e_2)} \ar[d]_p & \PP^2 \times \PP^2 \\
\overline{M}_{0,2}(\PP^2, 1) \ar[ur]_{(e_1, e_2)} .
} \end{equation*}
Hence
\begin{equation*}
\GW{\gamma_2, \gamma_2}{(1,0)} = \int (e_1, e_2)^\ast(h_1^2 h_2^2) = \GW{h^2, h^2}{1}^{\PP^2} = 1 .
\end{equation*}
The is the only invariant of this type that we will need.

\paragraph{Comments on the remaining invariants}
\label{comments}

The other invariants of the form $\GW{\phi_1, \phi_2}{(1,g)}$ are omitted here because, as we will see in Proposition~\ref{prop:gen-invs},  they can be deduced via the WDVV equations from the invariants we have already calculated.

It is less tedious in practice, however, to calculate these invariants directly.  When $g \leq 0$, they can be computed by translating them into questions about lines in $\PP^2$, as was done for the invariants above.

When $g > 0$, the moduli spaces $\overline{M}_n(1, g)$ have excess dimension and computing the invariants requires a virtual class calculation.  If the cycles $\phi_1$ and $\phi_2$ are chosen appropriately, then  $e^{-1}(\phi_1 \times \phi_2) \subset \overline{M}_n(1,g)$ is contained in the locus of comb curves with exactly $2$ orbifold points on the handle.  Thus, when $g > 0$, $\GW{\phi_1, \phi_2}{(1,g)}$ may be calculated as a sum over the partitions ${\bf h}$ of $[2g + 2]$ into two subsets, each containing an odd number of elements, of contributions from the $U_n(1, {\bf h})$.  By Theorem~\ref{thm:mod:1}, this contribution is precisely $(-\frac{1}{4})^g$ whenever it is nonzero, and one need only count the number of contributing partitions.  One obtains, for $g > 0$,
\begin{equation*}
\GW{\phi_1, \phi_2}{(1,g)} = \begin{cases} (-1)^g \GW{\phi_1, \phi_2}{(1,0)} & \text{at least one $\phi_i$ is untwisted} \\ (-1)^g \frac{1}{2} \GW{\phi_1, \phi_2}{(1,0)} & \text{both $\phi_i$ are twisted}.
\end{cases}
\end{equation*}

\subsection{Curves of a given degree through ordinary points}
\label{n-point-deg-d}

As in Section~\ref{comb-curves} we take $[2g + 2] = \set{1, 2, \ldots, 2g + 2}$.  Let $H$ be the set of all partitions of $[2g + 2]$ such that every part has an odd number of elements.  For each ${\bf h} \in H$, define $g({\bf h})$ to be the number $g'$ such that $2g' + 2$ is the number of subsets in the partition ${\bf h}$.

Let $P_1, \ldots, P_{3d + 1}$ be generic points in $\PP^2$.  We define 
\begin{equation*} \begin{split}
  B_i & = [ (P_i \times \PP^2 \cup \PP^2 \times P_i) / S_2 ] \subset [\Sym^2 \PP^2]  \\
  B & = B_1 \times \cdots \times B_{3d + 1} \subset [\Sym^2 \PP^2]^{3d + 1}
  \end{split}
\end{equation*}
Let $E(d,g)$ be the number of hyperelliptic curves of genus $g$ in $\PP^2$ passing through the points $P_1, \ldots, P_{3d + 1}$.

\begin{theorem} \label{thm:2}
  The following relationship between the Gromov--Witten invariants and enumerative invariants holds.
  \begin{equation}\label{eqn:gw}
    \left< B \right>_{(d,g)}
    = \sum_{{\bf h} \in H} \left( - \frac{1}{4} \right)^{g - g({\bf h})} (2 g ({\bf h}) + 2)! \: E(d, g({\bf h}))
  \end{equation}
\end{theorem}

A more explicit version of Theorem \ref{thm:2} is
\begin{corollary}
  \begin{equation*}
    \left< B \right>_{(d,g)}
    = \sum_{g' \geq 0}  \sum_{\substack{b_1 + 2 b_2 + 3 b_3 + \cdots = g - g' \\ b_0 + b_1 + b_2 + b_3 + \cdots = 2 g' + 2}} \left( - \frac{1}{4} \right)^{g - g'} \frac{(2g + 2)!}{1!^{b_0} 3!^{b_1} 5!^{b_2} \cdots} \frac{(2 g' + 2)!}{b_0! b_1! b_2! b_3! \cdots} E(d, g').
  \end{equation*}
\end{corollary}
\begin{proof}
  We reorganize the sum over all partitions in the statement of the theorem as the sum, first over the number of parts in the partition, then over all partition with types having that many parts.  If $(1^{b_0} 3^{b_1} 5^{b_2} \cdots)$ is a partition type with $g'$ parts then $\sum b_i = g'$ and $\sum (2i + 1) b_i = 2 g + 2$.  The second condition can be rewritten (using the first):
  \begin{equation*}
        2 g + 2 = \sum (2i + 1) b_i = 2 \sum i b_i + \sum b_i = 2 \left( \sum i b_i + g' + 1 \right)
  \end{equation*}
  so it is equivalent to the condition $\sum i b_i = g - g'$.
  This explains the indexing of the sum.  It remains to check that the number of partitions of $[2g + 2]$ with type $(1^{b_0} 3^{b_1} 5^{b_2} \cdots)$ is
  \begin{equation*}
    \left( \frac{(2g + 2)!}{1^{b_0} 3^{b_1} 5^{b_2} \cdots} \right) \left( \frac{1}{b_0! b_1! b_2! b_3! \cdots} \right) .
  \end{equation*}
  We may recognize the factor on the left as the number of partitions of $[2g + 2]$ into parts of odd orders, together with an ordering of the parts.  The factor on the right is simply the reciprocal of the number of ways of reordering the parts.
\end{proof}

Let $e : \overline{M}_n(d,g) \rightarrow [\Sym^2 \PP^2]^n$ be the map that evaluates the orbifold stable map at the ordinary marked points.  Then the Gromov--Witten invariant \ref{eqn:gw} is the virtual degree of
\begin{equation*}
  \overline{\Gamma}(d,g) = e^{-1}(B) \subset \overline{M}_{3d + 1}(d,g) .
\end{equation*}

\begin{proposition} 
  Every curve in $\overline{\Gamma}(d,g)$ is a comb curve.
\end{proposition}
\begin{proof}
  By \ref{dim-ev-image}, the image of the evaluation map $\overline{M}(d,g) \rightarrow [\Sym^2 \PP^2]^{3d + 1}$ has dimension $6d + 2$.  Let $Z$ be the image of $\overline{M}(d, g) - U(d,g)$ in $[\Sym^2 \PP^2]^{3d + 1}$.  By \ref{dim-ev-image}, $Z$ has image strictly less than $6d + 2$.  Since $[\Sym^2 \PP^2] - \Delta$ is homogeneous, and each $B_i$ has codimension $2$, the expected dimension of the intersection of $Z$ with $B \cap ([\Sym^2 \PP^2] - \Delta)$ is less than $6d + 2 - 2(3d + 1) = 0$, hence is empty when the $B_i$ are generic.  Therefore any intersection between $Z$ and $B$ must occur inside $\Delta^{3d + 1}$.  But $\Delta$ is also homogeneous and $B_i \cap \Delta$ has codimension $2$, so the same argument applies to show that $B \cap Z$ is empty when the $B_i$ are generic.  Thus the pre-image of $Z$ under the evaluation map is contained in $U(d,g)$.
\end{proof}

Define $\Gamma(d,g) = \overline{\Gamma}(d,g) \cap M_{3d + 1}(d,g)$ to be the substack of $\overline{\Gamma}$ parameterizing {\em smooth} orbifold curves interpolating the $B_i$.  By the proposition, $\Gamma(d,g)$ is contained in $U_{3d + 1}(d,g) \cap M_{3d + 1}(d,g) = M^\circ_{3d + 1}(d,g)$, the locus of curves smooth curves which meet the diagonal transversally.  By \ref{eval-is-smooth}, $M^\circ_{3d + 1}(d,g)$ is smooth.

\begin{proposition} \label{gamma-no-bar}
  The stack $\Gamma(d,g)$ is a disjoint union of $(2g + 2)! E(d,g)$ reduced points with only trivial automorphisms.
\end{proposition}
\begin{proof}
  We show first that $\Gamma(d,g)$ is a finite set of reduced points.  For this, note that $[\Sym^2 \PP^2]^{3d + 1}$ has finitely many orbits under the action of $\PGL{3}^{3d + 1}$.  Since $B$ is smooth and meets the orbit stratification in the expected dimension, $e^{-1}(B)$ will be smooth of the expected dimension when $B$ is chosen generically (by Kleiman--Bertini~\cite{Kl}; see also~\cite{G}, Lemma 2.5).  This implies that $\Gamma(d,g)$ has dimension zero (since $B$ has codimension $6d + 2$) and hence is a finite set of reduced points.  

  To prove the statement about automorphisms, first note that since $\Gamma(d,g) \subset M_{3d + 1}^\circ(d,g)$, if $(C, f) \in \Gamma(d,g)$ then $f$ does not factor through $\Delta$.  Thus $f$ does not carry the generic point of $C$ into $\Delta$, so $(C, f)$ can have an automorphism only if $f$ is a multiple cover of some curve $C'$ which is generically embedded in $[\Sym^2 \PP^2]$.  But $C$ has genus zero, so $C'$ has genus zero also, and therefore gives a point in $M^\circ_{3d' + 1}(d',g')$ for some $d' < d$ and $g' < g$.  On the locus in $M^\circ_{3d + 1}(d,g)$ of multiple covers factoring through such $C'$, the evaluation map to $[\Sym^2 \PP^2]^{3d + 1}$ will factor through $M^\circ_{3d' + 1}(d', g')$ and therefore have dimension at most $6d' + 2 < 6d + 2$.  Since there are only finitely many possibilities for $d'$ and $g'$, this means there is a closed subset of $[\Sym^2 \PP^2]$ of codimension greater than $6d + 2$.  Since $B$ meets the diagonal in the expected codimension, it follows that $e^{-1}(B)$ will not meet the locus of multiple covers if $B$ is generic.  Thus the points of $e^{-1}(B)$ will have no automorphisms.

  Now we argue that the number of these points is $(2g + 2)! E(d,g)$.  The moduil space $M_n^\circ(d,g)$ may also be viewed as the moduli space of smooth hyperelliptic curves in $\PP^2$ with $n$ marked pairs of hyperelliptically conjugate points and an ordering on the $2g + 2$ hyperelliptic branch points.  The evaluation map $e : M_n^\circ(d,g) \rightarrow [\Sym^2 \PP^2]^n$ is evaluation at the marked conjugate pairs.  Thus $\Gamma(d,g)$ is exactly the moduli space of hyperelliptic curves meeting the points $(P_1, \ldots, P_{3d + 1})$ with an ordering on the branch points.  The number of such curves is $   (2g + 2)! E(d,g)$.
\end{proof}
    
To prove the theorem it remains to understand the contributions of the compactification to \ref{eqn:gw}.  Let us write $\Gamma(d, {\bf h}) = \overline{\Gamma}(d,g) \cap U_{3d + 1}(d, {\bf h})$ (Section~\ref{comb-curves}).  Since $\overline{\Gamma}(d, {\bf h})$ is contained in $U_{3d + 1}(d, g)$, Corollary~\ref{comb-curves-disj-union} implies that
\begin{equation*}
  \overline{\Gamma}(d,g) = \coprod_{{\bf h} \in H} \Gamma(d, {\bf h}) .
\end{equation*}
By Theorem \ref{thm:mod:1}, the map $\Gamma(d, {\bf h}) \rightarrow \Gamma(d, g({\bf h}))$ has virtual degree $\left(-\frac{1}{4}\right)^{g - g({\bf h})}$.  Combining this with \ref{gamma-no-bar}, we obtain
\begin{equation*}
  \vdeg \Gamma(d, {\bf h}) = \left(-\frac{1}{4}\right)^{g - g({\bf h})} (2 g({\bf h}) + 2)! \: E(d, g({\bf h}) 
\end{equation*}
and summing over ${\bf h}$ gives
\begin{equation*}
  \left< B \right>_{(d, g)} = \sum_{{\bf h} \in H} \vdeg \Gamma(d, {\bf h}) = \sum_{{\bf h} \in H} \left(- \frac{1}{4}\right)^{g - g({\bf h})} (2 g({\bf h}) + 2)! \: E(d, g({\bf h})) .
\end{equation*}
which completes the proof of Theorem \ref{thm:2}.

%% file: gw-2008.06.15.tex
In this section, we will identify and calculate a collection of initial data that determine all of the genus zero Gromov--Witten invariants of $[\Sym^2 \PP^2]$.

\begin{proposition} \label{prop:gen-invs}
All of the orbifold Gromov--Witten invariants of $[\Sym^2 \PP^2]$ are determined from the structure of the orbifold Chow ring, together with the invariants
\begin{gather*}
\GW{\alpha \gamma, \gamma, \ldots, \gamma}{0} \\
\GW{\gamma^2, \gamma, \ldots, \gamma}{0}
\end{gather*}
and the $2$-point invariants, $\GW{\phi_1, \phi_2}{1}$, by means of the WDVV equations, the unit and divisor axioms, the dimension axiom, and linearity.
\end{proposition}

We recall that the unit axiom gives
\begin{equation*}
\GW{1, \phi}{d} = 0
\end{equation*}
if $\phi$ involves at least $3$ insertions or $d > 0$.  The divisor axiom gives
\begin{equation*}
\GW{\alpha, \phi}{d} = d \GW{\phi}{d}
\end{equation*}
under the same hypotheses.  By the dimension axiom, we mean the property that
\begin{equation*}
\GW{\phi_1, \ldots, \phi_n}{d} = 0
\end{equation*}
unless $\sum \deg(\phi_i) = 3d + 1 + n$.  The divisor and unit axioms are proven in~\cite{AGV}.  The dimension axiom is easy to deduce from the fact that the virtual fundamental class has degree $3d + 1 + n_0$, where $n_0$ is the number of ordinary marked points, and the orbifold degree of $\phi_i$ is its usual degree plus the age of the corresponding component of $\overline{I} [\Sym^2 \PP^2]$.

The proposition can be deduced from Proposition~\ref{prop:gen-q-invs}, so we will defer the proof.

\subsection{The orbifold Chow ring}
\label{chow-ring}
\input{chow-ring-2008.06.15.tex}

\subsection{The remaining invariants}

We have
\begin{equation*} \begin{split}
\GW{\alpha \gamma, \gamma, \ldots, \gamma}{0} & = 2 \GW{\gamma_1}{(0,g)} =  \frac{ (-1)^g \: (2^{2g} - 1) \: 6 \: B_{2g}}{2g} \\
\GW{\gamma^2, \gamma, \ldots, \gamma}{d} & = \GW{\gamma^2}{(0,g)} = 0 .
\end{split} \end{equation*}
These were both calculated in Section~\ref{deg-0}.

A $2$-point invariant, in order to be nonzero, must have an even number of orbifold points.  The only such invariants that are not zero for dimension reasons are
\begin{equation*} \begin{split}
\GW{\alpha^4, \alpha^2}{1} & = \GW{\alpha^4, \alpha^2}{(1,-1)} = 6 \\
\GW{\alpha^4, \gamma^2}{1} & = \GW{\alpha^4, \alpha^2 - \beta}{(1,-1)} = 6  \\
\GW{\alpha^3, \alpha^3}{1} & = \GW{\alpha^3, \alpha^3}{(1,-1)} = 9 \\
\GW{\alpha^2 \gamma, \alpha^2 \gamma}{1} & = 16 \GW{\gamma_2, \gamma_2}{(1,0)} = 16 
\end{split} \end{equation*}
which were all calculated in Section~\ref{2-point}.  To facilitate comparison in Section~\ref{crc-2-points} we also include
\begin{equation*} \begin{split}
\GW{\alpha^4, \alpha \gamma}{1} & = 0 \\
\GW{\alpha^3, \alpha^2 \gamma}{1} & = 0 .
\end{split} \end{equation*}

%% file: chow-ring-2008.06.15.tex
We calculated the group structure of $A^\ast_{\rm orb}([\Sym^2 \PP^2])$ in \ref{group-structure}, so we only need to understand the product.  We begin by recalling the definition of the orbifold product and orbifold Poincar\'e pairing.

Let $\phi_1, \phi_2$ be classes in $A^\ast_{\rm orb}(X)$ and let $M_{0,3}(X, 0)$ be the moduli space of $3$-pointed, degree zero, genus zero orbifold stable maps to $X$ (with arbitrary stack structure at the marked points).  The definition of the orbifold product in \cite{AGV} is
\begin{equation*}
  \phi_1 . \phi_2 = r (e_3)_\ast ( e_1^\ast(\phi_1) e_2^\ast(\phi_2) )
\end{equation*}
where $r$ is the order of the automorphism group at the third marked point.  This definition is made exactly so that
\begin{equation} \label{pairing-compatibility}
  ( \phi_1 . \phi_2 , \phi_3 )_{\rm orb} = ( \phi_1, \phi_2 . \phi_3 )_{\rm orb} = \left< \phi_1, \phi_2, \phi_3 \right>_0
\end{equation}
where
\begin{equation*}
  (\phi_1, \phi_2)_{\rm orb} = \left< \phi_1, \phi_2, 1 \right>_0
\end{equation*}
is the orbifold Poincar\'e pairing.

When the orbifold Chow ring of $X$ satisfies Poincar\'e duality (as is the case when $X = [\Sym^2 \PP^2]$, Equation~\eqref{pairing-compatibility} implies that the degree zero invariants determine the product on $A^\ast_{\rm orb}(X)$ by means of the orbifold Poincar\'e pairing.  That is, $\phi_1 . \phi_2$ is the unique class in $A^\ast_{\rm orb}(X)$ such that
\begin{equation*}
  \GW{\phi_1 . \phi_2, \phi_3, 1}{0} = \GW{ \phi_1, \phi_2, \phi_3}{0} .
\end{equation*}
Conversely, the $3$-point, degree zero invariants can be extracted from the structure of the orbifold Chow ring, since $\GW{\phi_1, \phi_2, \phi_3}{0} = \int \phi_1 \phi_2 \phi_3$, where $\int$ is the $\QQ$-linear function
\begin{equation*}
\int : A^\ast_{\rm orb}([\Sym^2 \PP^2]) \rightarrow \QQ
\end{equation*}
taking the value zero on the untwisted sector and restricting to the usual integration map on $A_\ast([\Sym^2 \PP^2])$.

\subsubsection{The orbifold Poincar\'e pairing on $[\Sym^2 \PP^2]$}

We specialize now to the orbifold Poincar\'e pairing of $X = [\Sym^2 \PP^2]$.  Recall that the rigidified inertia stack is $\overline{\Omega}_0 \amalg \overline{\Omega}_1$, where $\overline{\Omega}_0 = [\Sym^2 \PP^2]$ and $\overline{\Omega}_1 = \PP^2$.  We note that if $\phi_1 \in A^\ast(\overline{\Omega}_0)$ and $\phi_2 \in A^\ast(\overline{\Omega}_1)$ then $(\phi_1, \phi_2)_{\rm orb} = 0$ because the Gromov--Witten invariant $\left< \phi_1, \phi_2, 1 \right>_0$ is evaluated on the substack of $\overline{M}_{0,3}([\Sym^2 \PP^2], 0)$ that parameterizes curves with only one orbifold point, and this substack is empty.

Thus the orbifold Poincar\'e pairing is the direct sum of a pairing on $A^\ast(\overline{\Omega}_0)$ and one on $A^\ast(\overline{\Omega}_1)$.  It is easy to show that if $\phi_1, \phi_2 \in A^\ast(\overline{\Omega}_0)$ then
\begin{equation*}
  (\phi_1, \phi_2)_{\rm orb} = \int_{\overline{\Omega}_0} \phi_1 \phi_2 = (\phi_1, \phi_2)
\end{equation*}
where $(\hphantom{-},\hphantom{-})$ is the usual Poincar\'e pairing on $\overline{\Omega}_0 \cong [\Sym^2 \PP^2]$.  For $\phi_1, \phi_2 \in A^\ast(\overline{\Omega}_1)$, consider the substack $M'$ of $\overline{M}_{0,3}([\Sym^2 \PP^2], 0)$ where the first two evaluation maps are in $\overline{\Omega}_1$.  This is where $(\phi_1, \phi_2)_{\rm orb}$ is computed.  The two maps
\begin{equation*}
  e_1, e_2 : M' \rightarrow \overline{\Omega}_1
\end{equation*}
coincide and make $M'$ into a $S_2$-gerbe over $\overline{\Omega}_1$.  Hence
\begin{equation*}
  (\phi_1, \phi_2)_{\rm orb} = \frac{1}{2} \int_{\overline{\Omega}_1} \phi_1 \phi_2 = \frac{1}{2} ( \phi_1, \phi_2 )
\end{equation*}
where $(\hphantom{-},\hphantom{-})$ is in this case the usual Poincar\'e pairing on $\overline{\Omega}_1 \cong \PP^2$.

We can now write down the matrix of the Poincar\'e pairing with respect to the basis given in Section~\ref{group-structure}.  It is below.

\begin{equation} \label{table} \begin{array}{c|ccccccccc}
    & 1 & \alpha & \alpha^2 & \beta & \alpha^3 & \alpha^4 & \gamma_0 & \gamma_1 & \gamma_2 \\
    \hline
    1 &  &  &  &  &  & 3 &  &  & \\
    \alpha & & & & & 3 & \\
    \alpha^2 & & & 3 & 1 \\
    \beta & & & 1 & 3 \\
    \alpha^3 & & 3 \\
    \alpha^4 & 3 \\
    \gamma_0 & & & & & & & & & \frac{1}{2} \\
    \gamma_1 & & & & & & & & \frac{1}{2} \\
    \gamma_2 & & & & & & & \frac{1}{2} 
  \end{array}
\end{equation}

\subsubsection{The product structure}

We need only compute the product $\phi_1 . \phi_2$ when at least one of the $\phi_i$ is in the twisted sector.  Suppose first that $\phi_2$ is and $\phi_1$ is not.  Then let $M$ be the locus in $\overline{M}_{0,3}([\Sym^2 \PP^2]$ where the second and third evaluation maps are in the untwisted sector.  Then $M \cong \Delta$ and 
\begin{equation*}
  \phi_1 . \phi_2 = 2 p_\ast(i^\ast(\phi_1) \phi_2)
\end{equation*}
where $i : \Delta \rightarrow [\Sym^2 \PP^2]$ is the inclusion and $p : \Delta \rightarrow \PP^2$ is the projection on the coarse moduli space.  In particular, we obtain
\begin{equation*} \begin{split}
    \alpha . \gamma_0 & = 2 \gamma_1 \\
    \alpha . \gamma_1 & = 2 \gamma_2
  \end{split}
\end{equation*}
so $A^\ast_{\rm orb}([\Sym^2 \PP^2])$ is generated by $\alpha$, $\gamma_0$, and $\beta$ as a ${\bf Q}$-algebra.  From now on, we will write $\gamma$ instead of $\gamma_0$ to remove some clutter from the notation.

If $\phi_1$ and $\phi_2$ both come from the twisted sector, then
\begin{equation*}
  \phi_1 . \phi_2 = i_\ast (p^\ast (\phi_1 \phi_2))
\end{equation*}
and therefore $\gamma^2$ is the class of the diagonal in $A^\ast([\Sym^2 \PP^2])$.  

Since the class of the diagonal may also be expressed as $\alpha^2 - \beta$, it now follows that $A^\ast_{\rm orb}([\Sym^2 \PP^2])$ is generated as a ring by $\alpha$ and $\gamma$.  Of course, there are algebraic relations among $\alpha$ and $\gamma$ in addition to the degree constraints.  For example,
\begin{gather}  \label{eqn:6}
3 \pi^\ast(\alpha \gamma^2) = 6(h_1^2 h_2 + h_1 h_2^2) = 2\pi^\ast(\alpha^3) \\
3 \alpha^2 \gamma = 12 \gamma_2 = 4 \gamma^3. \notag
\end{gather}
where $\pi : \PP^2 \times \PP^2 \rightarrow [\Sym^2 \PP^2]$ is the canonical projection in the first line.  We therefore have relations,
\begin{equation} \label{relations} \begin{split}
R_1 & = 2 \alpha^3 - 3 \alpha \gamma^2 \\
R_2 & = 3 \alpha^2 \gamma - 4 \gamma^3.
\end{split} \end{equation}

\begin{proposition}
The relations $R_1$ and $R_2$ generate all of the relations in $A^\ast_{\rm orb}([\Sym^2 {\bf P}^2])$ between~$\alpha$ and~$\gamma$.
\end{proposition}
\begin{proof}
Let $B = {\bf Q}[\alpha, \gamma] / (R_1, R_2)$ and let $A = A^\ast_{\rm orb}([\Sym^2 {\bf P}^2])$.  It will be enough to show that the dimensions of each of the graded pieces coincide.  Let $A_n$ and $B_n$ be the $n$-th graded pieces of $A$ and $B$, respectively.

Since the relations are only in degree $3$, we only have to check this in degrees $3$ and higher.  In degree $3$, there are two independent relations, so $\displaystyle \dim B_3 = 4 - 2 = 2$, which again coincides with $A$.  The element $\alpha^4$ spans $B_4$ so $\dim B_4 = \dim A_4 = 1$.  

It remains to show that $B_5 = 0$.  We need only check that $\alpha^5$ and $\alpha^4 \gamma$ are  $0$ in $B$.  Indeed, $\alpha R_1 + \gamma R_2 = 0$ gives $\gamma^4 = \frac{1}{2} \alpha^4$, but then
\begin{gather*}
\frac{1}{2} \alpha^4 \gamma = \gamma^5 = \frac{3}{4} \alpha^2 \gamma^3 = \frac{9}{16} \alpha^4 \gamma \qquad \text{and} \\
2 \alpha \gamma^4 = \alpha^5 = \frac{3}{2} \alpha^3 \gamma^2 = \frac{9}{4} \alpha \gamma^4
\end{gather*}
in $B$.  These imply that $\alpha^5 = \alpha^4 \gamma = 0$ in $B_5$, hence $B_n = 0$ for $n \geq 5$ and $B \rightarrow A$ is an isomorphism.
\end{proof}

%% file: crc-2008.06.15.tex
Ruan's crepant resolution conjecture predicts that any two crepant resolutions of the same singular space should have equivalent Gromov--Witten theories after an appropriate change of variables.  The orbifold $[\Sym^2 \PP^2]$ may be viewed as a crepant resolution of its coarse moduli space, which is the scheme $\Sym^2 \PP^2$ and has an $A_1$ singularity along the diagonal.  Any $A_1$ surface singularity admits a crepant resolution by blowing up the singularity.  In the case of $\Sym^2 \PP^2$, this produces the Hilbert scheme, $\Hilb_2 \PP^2$, whose genus zero Gromov--Witten invariants were computed by Graber~\cite{G}.  In this section we verify that the relationship between the Gromov--Witten invariants of $[\Sym^2 \PP^2]$ and $\Hilb_2 \PP^2$ predicted by the crepant resolution conjecture is correct.

The crepant resolution conjecture for orbifolds was first formulated by Ruan \cite{R} for the degree zero, genus zero Gromov--Witten invariants.  An observation of Perroni~\cite{P} indicated that it needed modification and a new statement of the conjecture was provided by Bryan and Graber \cite{BG}.  Their statement also extended the conjecture to all degrees, but was determined by Coates, Corti, Iritani, and Tseng \cite{CCIT06} not to be the correct formulation in the absensce of the {\em hard Lefschetz condition}.  Coates and Ruan give in \cite{CoR} an updated version of the conjecture which makes use of Givental's Lagrangian cone formalism and applies in all genera.  In the presence of the hard Lefschetz condition, it reduces to the statement of Bryan and Graber.

An orbifold $X$ satisfies the hard Lefschetz condition if the automorphism of $\overline{I} X$ that inverts the band preserves the age.  This is trivially verified when the stabilizer groups have order $2$, which is the case for $[\Sym^2 \PP^2]$.  Thus the Coates--Ruan version of the crepant resolution conjecture specializes to the Bryan--Graber version in our example.  We will specialize further in our statement to the case of $A_1$ singularities to simplify the exposition.

\newcommand{\CRCGW}[3]{\langle {#1} \rangle^\ast_{#2}({#3})}
\begin{conjecture}[Ruan \cite{R}, Bryan--Graber \cite{BG}, Coates--Ruan \cite{CoR}]
  Let $X$ be an orbifold all of whose stabilizer groups all have order~$2$.  Let $Z$ be the coarse moduli space of $X$ and assume that $Z$ has a crepant resolution,~$Y$.
  \begin{equation*} \xymatrix{
      X \ar[dr]_q & & Y \ar[dl]^p \\
      & Z
    }
  \end{equation*}
  Let $E$ be the exceptional divisor of $Y \rightarrow Z$.  Identify $H_2(X, \ZZ)$ with the subgroup of $\beta \in H_2(Y, \ZZ)$ such that $\beta . E = 0$.  Then
  \renewcommand{\labelenumi}{(\alph{enumi})}
  \begin{enumerate}
  \item  There is an isomorphism of graded vector spaces
    \begin{equation*}
      L : A^\ast_{\rm orb}(X) \rightarrow A^\ast(Y)
    \end{equation*}
    extending the homomorphism $A^\ast(X) \cong A^\ast(Z) \xrightarrow{p^\ast} A^\ast(Y)$.

  \item For any $\phi \in A^\ast_{\rm orb}(X^n)$ and $\beta \in H_2(X, \ZZ)$, the function
    \begin{equation*}
      \CRCGW{\phi}{\beta}{q} = \sum_{a \in \ZZ} \left< L(\phi) \right>^Y_{\beta + a E} q^a
    \end{equation*}
    is meromorphic near $q = 0$ and has analytic continuation to $q = -1$.   

  \item If $\phi \in A^\ast_{\rm orb}(X^n)$, then 
    \begin{equation*}
      \left< \phi \right>^X_\beta = \CRCGW{\phi}{\beta}{-1}.
    \end{equation*}
  \end{enumerate}
\end{conjecture}

We will prove this conjecture in the case $X = [\Sym^2 \PP^2]$ and $Y = \Hilb_2 \PP^2$ by reducing it to a small number of explicit checks using the WDVV equations.  Since the families of multilinear functions
$\GW{\hphantom{-}}{d}$ 
and 
$\CRCGW{{\hphantom{-}}}{d}{q}$
both satisfy the WDVV equations --- and since the WDVV equations for $\CRCGW{\hphantom{-}}{d}{q}$ reduce to valid equations for $\GW{\hphantom{-}}{d}$ --- it is sufficient to check the conjecture on any collection of invariants that determine all others by means of the WDVV equations {\em for $\CRCGW{\hphantom{-}}{d}{q}$}.

\subsection{The Hilbert scheme}

The results of this section are copied from~\cite{G}.

We view $H = \Hilb_2 \PP^2$ as a $\PP^2$-bundle over $G = \Grass(2, 3)$, the variety parameterizing lines in $\PP^2$, the projection being the map which sends a length-$2$ subscheme of $\PP^2$ to the unique line containing it.  The Chow ring of $\Hilb_2 \PP^2$ is generated by the Chern classes $T_1 = c_1({\cal O}_{G}(1))$ and $T_2 = c_1({\cal O}_{H / G}(1))$.  It is given by the relations,
\begin{equation} \label{presentation}
A^\ast(\Hilb_2 \PP^2) = {\bf Q}[T_1, T_2] / (T_1^3, T_2^3 - 3 T_1 T_2^2 - 3 T_1^2 T_2).
\end{equation}

 Let $B_1 \subset \Hilb_2 \PP^2$ be the locus of non-reduced length $2$ subschemes supported at that point in $\PP^2$.  Let $B_2 \subset \Hilb_2 \PP^2$ be the locus of length $2$ subschemes of $\PP^2$ that contain a fixed point and are contained in a fixed line.  These are curves in $\Hilb_2 \PP^2$ and
\begin{equation*}
  T_i . B_j = \delta_{ij} .
\end{equation*}

Let $E$ be the exceptional locus of the blow-up $\Hilb_2 \PP^2 \rightarrow \Sym^2 \PP^2$.  It parameterizes non-reduced length $2$ subschemes of $\PP^2$ and it is linearly equivalent to
\begin{equation*}
  (E . B_1) T_1 + (E . B_2) T_2 = -2 T_1 + 2 T_2 .
\end{equation*}
To see that $E . B_1 = -2$, note that $B_1 \subset \Delta$ and $\Hilb_2 \PP^2$ is a resolution of an $A_1$ singularity.  For $E . B_2 = 2$, note that all length $2$ subschemes in a fixed line are parameterized by the coefficients $(a, b, c) \in \PP^2$ of a degree $2$ polynomial in $2$ variables, $a x^2 + b xy + c y^2$.  The intersection with $\Delta$ is the vanishing of the discriminant, which has degree $2$.

We identify $H_2( \Sym^2 \PP^2, \ZZ)$ with the subgroup of $H_2(\Hilb_2 \PP^2, \ZZ)$ having zero intersection with $E$.  Since $E = 2 (T_2 - T_1)$, this is
\begin{equation*}
  H_2( \Sym^2 \PP^2, \ZZ) = \ZZ . (B_1 + B_2) .
\end{equation*}

The following table displays all of the Gromov--Witten invariants of $\Hilb_2 \PP^2$ that we will need to verify the crepant resolution conjecture.  It is given in \cite{G}, Section 4.1, following Theorem 4.2, in a different basis.
\begin{equation} \label{eqn:12} \begin{array}{|c|r|r|r|}
\hline 
\beta = & B_2 & B_1 + B_2 & 2 B_1 + B_2 \\
\hline 
\GW{T_2^2, T_2^4}{\beta} & 3 & 12 & 3 \\
\GW{T_2^3, T_2^3}{\beta} & 9 & 27 & 9 \\
\GW{(T_2 - T_1)^2, T_2^4}{\beta} & 3 & -9 & -6 \\
\GW{T_2^2 (T_2 - T_1), T_2^2 (T_2 - T_1)}{\beta} & 4 & -8 & 4 \\
\GW{T_2^4, T_2(T_2 - T_1)}{\beta} & 3 & 0 & -3 \\
\GW{T_2^3, T_2^3(T_2 - T_1)}{\beta} & \frac{1}{2} & 0 & -\frac{1}{2} \\
\hline 
\end{array} \end{equation}
These invariants all vanish for $\beta = a B_1 + B_2$ with $a > 2$.

We will also use the invariants,
\begin{equation} \label{eqn:13} \begin{split}
\GW{ T_2^2 - T_1 T_2}{a B_1} & = \frac{-6}{a^2} \\
\GW{ (T_2 - T_1)^2}{a B_1} & = \frac{-9}{a^2} \qquad \text{for $a \geq 1$}
\end{split}
\end{equation}

\subsection{Verification of the conjecture}

\paragraph{The Chow ring isomorphism}
\label{chow-isom}

For $\phi_1, \phi_2, \phi_3 \in A_\ast(\overline{I} [\Sym^2 \PP^2])$, we have
\begin{equation*}
\GW{\phi_1, \phi_2, \phi_3}{0} = \int \phi_1 \phi_2 \phi_3
\end{equation*}
where the product is the orbifold product, and the integral symbol stands for the degree map on the untwisted sector and zero elsewhere.  Thus, the degree zero invariants are encoded in the structure of the orbifold Chow ring and the above integration map.  To check the crepant resolution conjecture for the $3$-point, degree~$0$ invariants, it will therefore be sufficient to check that $L$ determines an isomorphism between the orbifold Chow ring of  $[\Sym^2 \PP^2]$ and the  quantum corrected Chow ring of $\Hilb_2 \PP^2$.

In fact, this has already been proved by Perroni~\cite{P}, who actually proved the corresponding assertion in general for $A_1$ singularities, so the explicit verification given below is nothing new.  However to check the Chow rings are isomorphic requires little beyond writing the definitions in this case, so we carry out the verification anyway.

By definition, the quantum corrected Chow ring of $\Sym^2 \PP^2$ has the multiplication, $\phi_1 . \phi_2 = \phi_1 \ast \phi_2 \rest{q = -1}$ where
\begin{equation} \label{quantum-corrected-product}
  \phi_1 \ast \phi_2 = \sum_{a = 0}^\infty \left< \phi_1, \phi_2, \ast \right>^{\Hilb_2 \PP^2}_{a B_1} q^a .
\end{equation}
Of course, one must check that the substitution $q = -1$ is defined.  From~\cite{G}, Section 4.3, we have the relations (the quantum parameter $q_2$ from~\cite{G} is set to $0$ here),
\begin{equation*} \begin{split} 
    T_1 \ast T_1 \ast T_1 - 9 f^2 T_1 \ast T_2 \ast T_2 + (9 f^2 - 2 f) T_2 \ast T_2 \ast T_2  = 0 \\
    (1 - 18 f) T_2 \ast T_2 \ast T_2 - 3(1 - 6f) T_1 \ast T_2 \ast T_2 + 6 T_1 \ast T_1 \ast T_2 = 0
\end{split} \end{equation*}
where $f = \frac{q}{1 - q}$.  Hence the relations in the quantum corrected Chow ring are generated by
\begin{equation*} \begin{split}
    S_1 & = T_1 . T_1 . T_1 - \frac{9}{4} T_1 . T_2 . T_2 + \frac{13}{4} T_2 . T_2 . T_2  \\
    S_2 & = 5 T_2 . T_2 . T_2 - 6 T_1 . T_2 . T_2 + 3 T_1 . T_1 . T_2
\end{split} \end{equation*}
It is now easy to check that
\begin{equation*} \begin{split}
L : A^\ast_{\rm orb}([\Sym^2 {\bf P}^2]) \tensor {\bf C} & \rightarrow A^\ast(\Hilb_2 {\bf P}^2) \tensor {\bf C} \\
\alpha & \mapsto T_2 \\
\gamma & \mapsto i(T_2 - T_1)
\end{split} \end{equation*}
is an isomorphism, since
\begin{equation*} \begin{split}
  L(R_1) & = S_2 \\ 
  L(R_2) & = 4 i (S_2 - S_1)
\end{split} \end{equation*}
(recall from \eqref{relations} that $R_1$ and $R_2$ are the relations in $A^\ast_{\rm orb}([\Sym^2 \PP^2])$).

We also check that $L$ preserves the integration map.  We have $\int_{[\Sym^2 \PP^2]} \alpha^4 = 3$ and $\int_{\Hilb_2 \PP^2} L(\alpha^4) = \int_{\Hilb_2 \PP^2} T_2^4$.  Recall that $T_2$ is the divisor in $\Hilb_2 \PP^2$ of length $2$ subschemes of $\PP^2$ incident to a fixed line; thus $T_2^4$ is the locus of length $2$ subschemes incident to $4$ fixed lines.  There are $3$ such, so the integration maps coincide.

From now on we will identify the $\QQ[i]$-vector spaces, $A^\ast_{\rm orb}([\Sym^2 \PP^2]) \tensor[\QQ] \QQ[i]$ and $A^\ast(\Hilb_2 \PP^2) \tensor[\QQ] \QQ[i]$ by means of $L$. We will thus speak of the product $\phi_1 \ast \phi_2 = L^{-1}(L(\phi_1) \ast L(\phi_2))$ for $\phi_1, \phi_2 \in A^\ast_{\rm orb}([\Sym^2 \PP^2]) \tensor[\QQ] \QQ[i]\left[[q]\right]$.  We will also write $\phi^{\ast n}$ to mean $\phi \ast \cdots \ast \phi$ ($n$ times).

\paragraph{The WDVV equations and divisor axiom}

Here we prove some properties of the invariants~$\CRCGW{\hphantom{-}}{d}{q}$ that we will use in a moment.

The divisor axiom and WDVV equations for $\CRCGW{\phi}{d}{q}$ follow from the corresponding properties of the Gromov--Witten invariants of $\Hilb_2 \PP^2$.  For the divisor axiom, we have
\begin{equation*}
\CRCGW{\alpha, \phi}{d}{q} = \sum_{a \in \ZZ} \CRCGW{T_2, L(\phi)}{(a + d) B_1 + d B_2}{q} = d \sum_{a \in \ZZ} \CRCGW{L(\phi)}{(a + d) B_1 + d B_2}{q} = d \CRCGW{\phi}{d}{q} .
\end{equation*}
Similarly, we have the unit axiom: if $d > 0$ or $\phi$ has at least $3$ insertions, then
\begin{equation*}
\CRCGW{1, \phi}{d}{q} = \sum_{a \in \ZZ} \CRCGW{1, L(\phi)}{(a + d) B_1 + d B_2}{q} = 0 .
\end{equation*}
For the WDVV equations, let $\set{\xi_i}$ be a homogeneous basis of $A_\ast(\overline{I} [\Sym^2 \PP^2])_{\QQ\left[[q]\right]}$ and let $\set{\tilde{\xi}_i}$ be the dual basis of $A_\ast(\overline{I} [\Sym^2 \PP^2])_{\QQ\left[[q]\right]}$ with respect to the pairing $\phi_1 \tensor \phi_2 \mapsto \CRCGW{\phi_1, \phi_2, 1}{0}{q}$ (note that this pairing is non-degenerate because when $q = 0$ it is the Poincar\'e pairing on $\Hilb_2 \PP^2$).  We have
\begin{multline*}
\sum_{\substack{S_1 \amalg S_2 = [n] \\ 1, 2 \in S_1 ; 3, 4 \in S_2 \\ d_1 + d_2 = d}} \sum_{i,j} \CRCGW{\phi_{S_1}, u_i}{d_1}{q} g_{ij} \CRCGW{\phi_{S_2}, u_j}{d_2}{q}  \\
= \sum_{a_1, a_2 \in \ZZ} \sum_{\substack{S_1 \amalg S_2 = [n] \\ 1, 2 \in S_1 ; 3, 4 \in S_2 \\ d_1 + d_2 = d}} \sum_{i,j} \GW{L(\phi_{S_1}), L(\xi_i)}{(a_1 + d_1) B_1 + d_1 B_2} \GW{L(\phi_{S_2}), L(\tilde{\xi}_j)}{(a_2 + d_2) B_1 + d_2 B_2} q^{a_1 + a_2} .
\end{multline*}
We may apply the WDVV equations on $\Hilb_2 \PP^2$ to the second line and then reverse the above equality to deduce the WDVV equations for the functions $\CRCGW{\phi}{d}{q}$.

Note that if the substitution $q = -1$ is legitimate, then the WDVV equations, the divisor axiom, and the unit axiom all reduce to the corresponding equations ans axioms for the Gromov--Witten invariants of $[\Sym^2 \PP^2]$.  Indeed, by Section~\ref{chow-isom} the $\xi_j$ reduce to a basis of $A_\ast(\overline{I} [\Sym^2 \PP^2])$ over $\QQ$ and the $\tilde{\xi}_j$ reduce to a dual basis since the pairing $\CRCGW{\phi_1, \phi_2,1}{0}{-1}$ reduces to the orbifold Poincar\'e pairing.

\paragraph{Reduction to degree $0$ and $2$-point invariants}

\begin{proposition} \label{prop:gen-q-invs}
The values of all of the $\CRCGW{\phi_1, \dots, \phi_n}{d}{q}$ are determined by means of the WDVV equations, the divisor, unit, and dimension axioms, and linearity from the invariants
\begin{gather}
\CRCGW{\alpha \ast \gamma, \gamma, \ldots, \gamma}{0}{q} \label{inv1} \\
\CRCGW{\gamma^{\ast 2}, \gamma, \ldots, \gamma}{0}{q} \label{inv2}
\end{gather}
and the invariants $\CRCGW{\phi_1, \phi_2}{1}{q}$.
\end{proposition}

\begin{proof}
 The WDVV equations and divisor axiom tell us that
\begin{equation*}
\CRCGW{\alpha \ast \phi_1, \phi_2, \phi_3, \ldots}{d}{q} - \CRCGW{\alpha \ast \phi_2, \phi_1, \phi_3, \ldots}{d}{q}
\end{equation*}
can be expressed as a polynomial in terms of {\em lower invariants}, i.e., invariants with a smaller number of insertions or of smaller degree.  Thus, any invariant with at least $3$ insertions is determined by the invariants with fewer insertions and the invariants,
\begin{equation} \label{an-inv}
\CRCGW{\alpha^{\ast n} \ast \phi_1, \phi_2, \phi_3, \ldots, \phi_n}{d}{q}
\end{equation}
where each $\phi_i$ is either $\gamma$ or $\gamma^{\ast 2}$.  Using the WDVV equations, we will now show by induction on $n$, $d$, and the number of appearances of $\gamma^{\ast 2}$ among the $\phi_i$ that all of these invariants can be obtained as values of polynomials in the invariants~\eqref{inv1} and~\eqref{inv2}.

In the invariant~\eqref{an-inv}, if one of the $\phi_i$ is $\gamma^{\ast 2}$ then we may assume, using WDVV, that it is $\phi_1$.  If in addition $n > 0$, then we know $\alpha \ast \gamma^{\ast 2}$ has degree $3$, hence is a linear combination of $\alpha^{\ast 3}$ and $\alpha^{\ast 2} \ast \gamma$, because $\alpha^3$ and $\alpha^2 \gamma$ span $A^\ast_{\rm orb}([\Sym^2 \PP^2]) = A^\ast_{\rm q}(\Hilb_2 \PP^2) \pmod{q + 1}$ by Section~\ref{chow-isom}.  This permits us to reduce the number of appearances of $\gamma^{\ast 2}$.

We may now assume that either $\phi_i = \gamma$ for all $i$, or that $n = 0$.  In the former case, the invariant is zero by the dimension axiom unless $n = 1$, in which case the invariant is
\begin{equation*}
\CRCGW{\alpha \ast \gamma, \gamma, \ldots, \gamma}{d}{q} .
\end{equation*}
Otherwise, we consider invariants
\begin{equation} \label{an-inv2}
\CRCGW{\gamma^{\ast 2}, \ldots, \gamma^{\ast 2}, \gamma, \ldots, \gamma}{d}{q} .
\end{equation}
Assume $\gamma^{\ast 2}$ appears at least three times here.  Then by the WDVV equations,
\begin{multline*} 
\CRCGW{\gamma \ast \gamma, \gamma^{\ast 2}, \gamma^{\ast 2}, \ldots}{d}{q} + \CRCGW{\gamma, \gamma, \gamma^{\ast 2} \ast \gamma^{\ast 2}, \ldots}{d}{q} \\
\equiv \CRCGW{\gamma \ast \gamma^{\ast 2}, \gamma, \gamma^{\ast 2}, \ldots}{d}{q} + \CRCGW{\gamma, \gamma^{\ast 2}, \gamma \ast \gamma^{\ast 2}, \ldots}{d}{q} \pmod{\text{lower invariants}}.
\end{multline*}
Now, $\gamma^{\ast 3}$ is a linear combination of $\alpha^{\ast 3}$ and $\alpha^{\ast 2} \ast \gamma$, again because $\alpha^3$ and $\alpha^2 \gamma$ span $A^3_{\rm orb}([\Sym^2 \PP^2])$.  This allows us to reduce the number of appearances of $\gamma^{\ast 2}$ by one on the right side.  On the left side, $\gamma^{\ast 4}$ is proportional to $\alpha^{\ast 4}$ (because $\alpha^4$ spans $A^4_{\rm orb}([\Sym^2 \PP^2])$).  Taken together these give an expression of the invariant~\eqref{an-inv2} as a polynomial combination of invariants with fewer appearances of $\gamma^2$.

Thus we see that every invariant can be expressed as a polynomial combination of the $2$-point invariants and the invariants
\begin{gather*}
\CRCGW{\alpha \ast \gamma, \gamma, \ldots, \gamma}{d}{q} \\
\CRCGW{\gamma^{\ast 2}, \gamma, \ldots, \gamma}{d}{q} \\
\CRCGW{\gamma^{\ast 2}, \gamma^{\ast 2}, \gamma, \ldots, \gamma}{d}{q} .
\end{gather*}
We can easily see by the dimension axiom that the last of these is zero and the first two will be zero unless $d = 0$.  This completes the proof.
\end{proof}

Now to prove the crepant resolution conjecture, it will be sufficient to show that each of the $2$-point invariants and the invariants~\eqref{inv1} and~\eqref{inv2} is a meromorphic function of $q$ at $q = 0$ admitting analytic continuation to $q = -1$, and that when the substitution $q = -1$ is made, that invariant takes the same value as the corresponding invariant of $[\Sym^2 \PP^2]$.

\paragraph{The $2$-point invariants}
\label{crc-2-points}

The invariants that do not vanish by the dimension axiom are calculated below using Table~\ref{eqn:12}.
\begin{equation*} \begin{split}
  \CRCGW{\alpha^4, \alpha^2}{1}{q} & = \sum_{a \in \ZZ} \GW{T_2^4, T_2^2}{(a + 1)B_1 + B_2} q^a = 3 q^{-1} + 12 + 3q \\
  \CRCGW{\alpha^4, \gamma^2}{1}{q} & = \sum_{a \in \ZZ} \GW{T_2^4, \left(i (T_2 - T_1) \right)^2}{(a + 1)B_1 + B_2} q^a = -3 q^{-1} + 9 + 6 q \\
  \CRCGW{\alpha^3, \alpha^3}{1}{q} & = \sum_{a \in \ZZ} \GW{T_2^3, T_2^3}{(a + 1)B_1 + B_2} = 9 q^{-1} + 27 + 9 q \\
  \CRCGW{\alpha^2 \gamma, \alpha^2 \gamma}{1}{q} & = \sum_{a \in \ZZ} \GW{i T_2^2 (T_2 - T_1), i T_2^2 (T_2 - T_1)}{(a + 1)B_1 + B_2} = - 4 q^{-1} + 8 - 4q \\
  \CRCGW{\alpha^4, \alpha \gamma}{1}{q} & = \sum_{a \in \ZZ} \GW{T_2^4, i T_2(T_2 - T_1)}{(a + 1)B_1 + B_2} q^a = 3 i q^{-1} - 3 i q \\
  \CRCGW{\alpha^3, \alpha^2 \gamma}{1}{q} & = \sum_{a \in \ZZ} \GW{T_2^3, i T_2^2(T_2 - T_1)}{(a + 1)B_1 + B_2} q^a = \frac{i}{2} q^{-1} - \frac{i}{2} q 
\end{split} \end{equation*}

Substituing $q = -1$ and comparing with the calculations in Section~\ref{2-point-deg-1} completes the check for $2$-point invariants.

\paragraph{The degree zero invariants}

Now we show the agreement of the invariants,
\begin{gather*}
\GW{\alpha \gamma}{(0,g)} = \CRCGW{\alpha \gamma, \gamma, \ldots, \gamma}{0}{-1} \\
\GW{\gamma^2}{(0,g)} = \CRCGW{\gamma^2, \gamma, \ldots, \gamma}{0}{-1} 
\end{gather*}
where in each of the invariants on the left, $2g + 2$ is the number of appearances of $\gamma$ and $\alpha \gamma$ on the right.

The left sides of these equations were computed in Section~\ref{deg-0}.  For the right side, we have
\begin{equation*} \begin{split}
\CRCGW{\alpha \gamma, \gamma, \ldots, \gamma}{0}{q} 
& = \sum_{a = 1}^\infty \GW{i T_2(T_2 - T_1), i(T_2 - T_1), \ldots, i(T_2 - T_1)}{(a,0)} q^a \\
& = \sum_{a = 1}^\infty (-1)^{g} a^{2g + 1} \GW{T_2 (T_2 - T_1)}{(a,0)} q^a \\
& = \sum_{a = 1}^\infty (-1)^{g + 1} a^{2g - 1} 6 q^a \\
& = (-1)^{g + 1} 6 \Li_{-(2g - 1)}(q) \\
\CRCGW{\gamma^2, \gamma, \ldots, \gamma}{0}{q} 
& = \sum_{a = 1}^\infty \GW{- (T_2 - T_1)^2, i(T_2 - T_1), \ldots, i(T_2 - T_1)}{(a,0)} q^a \\
& = \sum_{a = 1}^\infty (-1)^g a^{2g + 2} \GW{(T_2 - T_1)^2}{(a,0)} q^a \\
& = \sum_{a = 1}^\infty (-1)^{g + 1} a^{2g} 9 q^a \\
& = (-1)^{g + 1} 9 \Li_{-2g}{q} .
\end{split} \end{equation*}
Evidently, these sums have analytic continuation to $q = -1$.  Recall now that $\Li_k(-1) = (2^{1 - k} - 1) \zeta(k)$.  Therefore,
\begin{equation*}\begin{split}
\CRCGW{\alpha \gamma, \gamma^{\tensor (2g + 1)}}{0}{-1} & = (-1)^{g + 1} 6 (2^{2g} - 1) \zeta(1 - 2g) = \frac{(-1)^g 6 (2^{2g} - 1) B_{2g}}{2g} \\
\CRCGW{\gamma^2, \gamma^{\tensor (2g+2)}}{0}{-1} & = (-1)^{g + 1} 9 (2^{2g + 1} - 1) \zeta(-2g) = 0 .
\end{split} \end{equation*}
These coincide with the calculations of Section~\ref{deg-0}.